\documentclass[11pt]{amsart}

\usepackage{amsfonts,newlfont,latexsym,pdfsync}
\usepackage{euscript}
\usepackage{amssymb,amscd}
\usepackage{amsmath}

      \theoremstyle{plain}
      \newtheorem{theorem}{Theorem}[section]
      \newtheorem{lemma}[theorem]{Lemma}
      \newtheorem{corollary}[theorem]{Corollary}
      \newtheorem{proposition}[theorem]{Proposition}
 \newtheorem{conjecture}{Conjecture}
   \newtheorem*{Theorem}{Main Theorem}

      \theoremstyle{definition}
        \newtheorem*{remark}{Remark}
        
        \newtheorem*{definition}{Definition}
        \newtheorem{problem}{Problem}
               \makeatletter
      \def\@setcopyright{}
      \def\serieslogo@{}
      \makeatother
\newcommand{\N}{\mathbb N}
\newcommand{\Z}{\mathbb Z}

\newcommand{\R}{\mathbb R}

\newcommand{\T}{\mathbb T}
\newcommand{\Zk}{\mathbb Z^k}
\newcommand{\rk}{\mathbb R^k}
\newcommand{\Tk}{\mathbb T^k}
\newcommand{\Tkp}{\mathbb T^{k+1}}

\newcommand{\Ce}{C^{1+\theta}}

\newcommand{\e}{\varepsilon}
\newcommand{\la}{\lambda}

\def \M{{\cal M}}

\def\dist{\text{dist}}

\def \a{\alpha}

\def \at{{\alpha (\mathbf t)}}
\def \as{{\alpha (\mathbf s)}}

\def \d{\delta}
\def \b{\beta}
\def \bt{{\beta (\mathbf t)}}
\def \bs{{\beta (\mathbf s)}}

\def \0{\mathbf 0}

\def \ev{\mathbf e}
\def \t{\mathbf t}
\def \s{\mathbf s}
\def \g{\mathbf g}
\def \gxt{{{\mathbf g} (x,\mathbf t)}}

\def \nor{\mathbf w}

\def \O{\cal{O}}
\def \on{\cal{N}}

\def \h{h}

\def \wl{W}
\def \w{{\cal W}}
\def \wbl{\tilde W}
\def \wb{\tilde {\cal W}}

\def \mw{\mu ^{\cal W}}

\def \de{D^E}

\def\proof{{\bf Proof.}\; }
\def \Re{\EuScript{R}}
\def \Rel{\EuScript{R}_\varepsilon ^l}
\def \P{\EuScript{P}}

\def \Qe{\EuScript{Q}}

\newcommand{\Id}{\operatorname{Id}}

\numberwithin{equation}{section}

\begin{document}

\author[Boris Kalinin, Anatole Katok, Federico Rodriguez Hertz]{Boris Kalinin$^{*})$, Anatole Katok $^{**})$, and\\ Federico Rodriguez Hertz$^{***})$}

\address {University of South Alabama, Mobile, AL}
\email{kalinin@jaguar1.usouthal.edu}

\address{The Pennsylvania State University, University Park, PA}
\email{katok\_a@math.psu.edu}

\address {Universidad de la Rep\'ublica, Montevideo, Uruguay}
\email{frhertz@fing.edu.uy}

\title[Nonuniform Measure Rigidity]{Nonuniform Measure Rigidity}
\dedicatory{Dedicated to the  memory of Bill Parry (1934--2006)}


\begin{abstract}
We consider an ergodic invariant measure $\mu$ for a smooth action $\a$
of  $\Zk,\,\,k\ge 2$, on a $(k+1)$-dimensional manifold  or for a locally free smooth action of $\rk,\,\,k\ge 2$ on a $(2k+1)$-dimensional manifold. We prove that
if $\mu$ is hyperbolic with the Lyapunov hyperplanes in general position
and if one element in $\Zk$ has positive entropy, then $\mu$ is absolutely
continuous. The main ingredient is absolute continuity of conditional measures on Lyapunov foliations which 
holds for a more general  class of smooth  actions of higher rank abelian groups. 

\end{abstract}


\date{\today}

\thanks{$^{*})$ Based on research supported by NSF grant DMS-0701292}
\thanks{$^{**})$ Based on research supported by NSF grants DMS-0505539 and DMS 0803880}
\thanks{$^{***})$ Based on research supported by PDT grants 54/18 and 63/204}

\maketitle


\section{Introduction}
In this paper we continue and significantly advance  the  line of development
started in \cite{KaK3} and \cite{KRH}. The general program is to show  that actions
of higher rank abelian groups, i.e  $\Zk\times\R^l,\,\,k+l\ge 2$,  by diffeomorphisms of compact
manifolds must preserve a geometric  structure, such as an absolutely continuous
invariant measure, under  global conditions of topological or dynamical  nature which
ensure both  infinitesimal   hyperbolic behavior and  sufficient global complexity of the
orbit structure.

In \cite{KaK3} and \cite{KRH} we considered $\Zk$ actions on the torus $\Tkp$, $k\ge 2$,
that induce on the first homology group  the action of a  maximal abelian subgroup
of  $SL(k+1,\Z)$ diagonalizable over $\R$. We say that such an  action has
{\em Cartan homotopy data}.\footnote{In the case of the torus it may look more natural
to speak about homology data, but we wanted to emphasize that  what mattered was
the homotopy types on individual elements; this notion can be generalized while
homological information in general is clearly insufficient.}
The central feature of that situation is existence of a semi-conjugacy $h$ between the
action, which we denote by $\a$, and the corresponding Cartan action $\a_0$  by
automorphisms of the torus, i.e. a  unique surjective continuous map $h : \Tkp \to \Tkp$
homotopic to identity  such that
$$
h \circ \a=\a_0 \circ h.
$$
This gives desired global complexity right away and allows to produce nonuniform
hyperbolicity (non-vanishing of the Lyapunov exponents) with little effort (see
\cite[Lemma 2.3]{KaK3}). Existence of  the semi-conjugacy  allows to use specific
properties of the  linear action $\a_0$ and reduces the proofs  to showing that the semi--conjugacy
is   absolutely continuous   and  bijective  on an invariant  set of positive Lebesque measure.\footnote{And in fact smooth in the sense of Whitney on smaller
               non-invariant sets of positive Lebesgue measure.}
Thus, this may be considered as a  version in the setting of global measure rigidity of the
{\em a priori} regularity method  developed for the study of local differentiable rigidity in
\cite{KS}  (see also an earlier paper \cite{KL}) and successfully  applied to the global
conjugacy problem on the torus in \cite{RH}.

In the present paper we consider an essentially different and more general situation.
We make no assumptions on the topology of the ambient manifold or the action under
consideration and instead
assume  directly that  the action preserves a measure  with non-vanishing Lyapunov
exponents  whose behavior is similar to that of the exponents for a Cartan action. 
Namely,  we consider a  $\Zk$, $k\ge 2$, action on a $(k+1)$-dimensional manifold 
or an $\rk ,\,\,k\ge 2$,  action on a $(2k+1)$-dimensional manifold  with an ergodic
invariant measure for which the kernels of the Lyapunov exponents are in general 
position (see the Definition below). Dynamical complexity is provided by the assumption 
that at least one element of the action has positive entropy. In fact our results for $\Zk$ actions are direct corollaries of those  for $\rk$ actions via suspension construction.

To formulate our results precisely recall that the Lyapunov characteristic exponents with
respect to an ergodic invariant measure for a smooth $\rk$ action are linear functionals
on $\R^k$. For a smooth $\Zk$ action they are linear functionals on $\Z^k$ which are 
extended to $\R^k$ by linearity. The kernels of these functionals are called the 
{\em Lyapunov hyperplanes}. A Lyapunov exponent is called {\em simple} if the
corresponding Lyapunov space is one-dimentional. See Section~\ref{section-prelim} 
for more details.

\begin{definition}
We will say that $m$ hyperplanes (containing $0$) in $\R^k$ are {\em in general position}
if  the dimension of the intersection of any $l$ of them is the minimal possible, i.e. is equal
to $\max\{ k-l, 0\}$.

We will say that the Lyapunov exponents of an ergodic invariant measure for a 
$\Zk$ action  are
{\em in general position} if they are all simple and nonzero, and if the Lyapunov
hyperplanes are distinct hyperplanes  in general position.

Similarly for an $\rk$ action the Lyapunov exponents of an ergodic invariant measure   are
{\em in general position} if the  zero exponent has multiplicity $k$  and 
the remaining exponents are all simple and nonzero, and if the Lyapunov
hyperplanes are distinct hyperplanes  in general position.

\end{definition}

\newpage

\begin{Theorem}\label{main}

\*

\begin{enumerate}
\item
Let $\mu$ be an ergodic invariant measure for a $\Ce,\, \theta >0$, action $\a$
of  $\Zk,\,\,k\ge 2$, on a $(k+1)$-dimensional manifold $M$. Suppose that the
Lyapunov exponents of $\mu$ are in general position and that at least one element
in $\Zk$ has positive entropy with respect to $\mu$.
Then $\mu$ is absolutely continuous.

\item
Let $\mu$ be an ergodic invariant measure for a  locally free $\Ce,\, \theta >0$, action $\a$
of   $\R^k,\,\,k\ge 2$, on a $2k+1$-dimensional manifold $M$. Suppose that 
Lyapunov exponents of $\mu$ are in general position and that at least one element
in $\R^k$ has positive entropy with respect to $\mu$.
Then $\mu$ is absolutely continuous.
\end{enumerate}
\end{Theorem}
\smallskip

As we already mentioned the  statement (1) is   a direct  corollary of (2)
applied to the  suspension of the  $\Zk$ action $\a$. We are not aware of any examples of $\rk$  actions satisfying assumptions of (2) other than  time changes of suspensions of $\Zk$  actions satisfying (1). 

Thus, what we prove is the first case of  existence of an absolutely continuous invariant
measure for actions of abelian groups whose orbits have codimension two or higher
which is derived  from general purely dynamical assumptions. Nothing of that sort takes
places in the classical dynamics for actions of orbit codimension two or higher.
\footnote{In our situation the  codimension of orbits  is at least three. When  codimension of orbits equals two  there is not  enough space for nontrivial behavior of  higher rank actions 
involving any kind of hyperbolicity, see \cite{K97}.} Only for codimension one actions (diffeomorphisms of the circle and fixed point free flows on the torus) of sufficient smoothness Diophantine condition on the rotation number (which is of dynamical nature)  guarantees
existence of a smooth invariant measure \cite{He, Y}.  One can point out  though that even
in those cases existence of topological conjugacy (for the circle) or orbit equivalence (for the
torus) to an algebraic system follows from the classical Denjoy theorem
(see e.g. \cite[Theorem 12.1.1]{KH}) and the
work goes into proving smoothness.  Thus  this falls under the general umbrella of {\em a priori}
 regularity methods, albeit   substantively very different from the hyperbolic situations, and should
 be   more appropriately compared with results of  \cite{KaK3} and \cite{KRH}.

In order to prove measure rigidity we develop principal elements of the
basic geometric approach  of \cite{KS3} in  this general non-uniform setting.
This has  been done partially already in \cite{KaK3} and we will rely on those
results  and constructions of that paper which do not depend on existence of
the semi-conjugacy.

The main technical problem which we face is showing  recurrence for  elements
within  the Lyapunov hyperplanes. For the actions on the torus
the semi-conjugacy  was used in a critical way.  One main innovation  here is a
construction of a particular time change which is smooth along the orbits of the
action but only measurable transversally  which ``straightens out'' the expansion
and contractions coefficients. This is somewhat similar to the ``synchronization''
time change for Anosov flows  introduced by  Bill Parry in \cite{Pa}. The main
technical difficulty lies in the fact that we need the  new action to possess certain
properties as if it were smooth. Section~\ref{TC} where this time change is defined
and its properties are studied  is the heart and the main technical part of the present
paper.

Results of this paper have been announced in \cite{KKRH-ERAMS}.



\section{Preliminaries}\label{section-prelim}

\subsection{Lyapunov exponents and suspension} \label{Lexponents}

In this section we briefly recall the definitions of Lyapunov characteristic exponents
and related notions for $\Zk$ and $\rk$ actions by measure preserving diffeomorphisms
of smooth manifolds. We refer to \cite[Sections 5.1 and 5.2]{KaK1} for more details
on general theory in the discrete case and to \cite{KaK3} for further development in
a more specialized setting. We will use those notions without special references.

Let $\a$ be a smooth $\Zk$ action on a manifold $M$ with an ergodic invariant
measure $\mu$. According to Multiplicative Ergodic Theorem for $\Zk$ actions
(see \cite{KaK1}) the Lyapunov decompositions for  individual elements of $\a$ have a
common refinement $TM =\bigoplus E_{\chi} $ called the {\em Lyapunov decomposition}
for $\a$. For each Lyapunov distribution $E_{\chi}$ the corresponding Lyapunov exponent,
viewed as a function of an element in $\Z^k$, is a linear functional $\chi : \Z^k \to \R$
which is called a {\em Lyapunov exponent} of $\alpha$. The Lyapunov exponents of
$\alpha$ are extended by linearity to functionals on $\rk$. The hyperplanes
$\,\ker \,\chi \subset \rk$ are called the {\em Lyapunov hyperplanes}  and the connected
components of $\,\rk \setminus \bigcup _{\chi} \ker \chi $ are called the {\em Weyl chambers}
of $\alpha$. The elements in the union of the Lyapunov hyperplanes are called 
{\em singular}, and the elements in the union of the Weyl chambers are called {\em regular}.
The corresponding notions for a smooth $\rk$ action are defined similarly (see 
Proposition \ref{MET} below for more details). We note that any $\rk$ action has 
$k$ identically zero Lyapunov exponents corresponding to the orbit directions.
These Lyapunov exponents are called {\em trivial} and the other ones are
called {\em nontrivial}. For the rest of the paper a Lyapunov exponent of an 
$\rk$ action will mean a nontrivial one.

One of the reasons for extending the Lyapunov exponents for a $\Zk$ action to $\rk$ is that
the Lyapunov hyperplanes may be irrational and hence ``invisible'' within $\Zk$. It is also
natural to construct an $\rk$ action for which the extensions of the exponents from $\Zk$
will provide the nontrivial exponents. This is given by the suspension construction which associates  to  a given $\Zk$ action on a manifold  $M$ an $\rk$ action on the suspension
manifold $S$, which is a bundle over $\Tk$ with fiber $M$. Namely, let $\Zk$ act on
$\rk \times M$  by $z\, (x,m) =(x-z,z \, m)$ and form the quotient space
\[ S = \rk \times M / \Zk.\]
Note that the action of $\rk$  on $\rk \times M$ by $x
\,(y,n)=(x+y,n)$ commutes with the $\Zk$-action and therefore
descends to $S$. This $\rk$-action is called the  {\em suspension} of the
$\Zk$-action.
There is a natural correspondence between the invariant measures, nontrivial Lyapunov
exponents, Lyapunov distributions, stable and unstable manifolds, etc. for the original
$\Zk$ action and its suspension. 

Since most of the arguments will be for the $\rk$ case, we summarize in the next 
proposition important properties of smooth $\rk$ actions given by the nonuniformly 
hyperbolic theory (see \cite{KaK1,BPbook}).
For a smooth $\rk$ action $\a$ on a manifold $M$ and an element $\t\in\rk$ we 
denote the corresponding diffeomorphism of $M$ by $\a(\t)$. Sometimes we will 
omit $\a$ and write, for example, $\t x$ in place of $\a(\t)x$ and $D \t$ in place of 
$D\a(\t)$ for the derivative of $\a(\t)x$.

 \begin{proposition} \label{MET}
Let $\a$  be a  locally free $\Ce,\, \theta >0$, action of $\R^k$ on a manifold $M$
preserving an ergodic invariant measure $\mu$.
There are linear functionals $\chi_i$, $i=1, \dots, l$, on $\rk$  and an $\alpha$-invariant
measurable splitting, called the {\em Lyapunov decomposition}, of the tangent bundle of
$M$
 $$TM = T\O \oplus \bigoplus_{i=1}^{l} E_{i}$$
over a set of full measure $\Re$, where $T\O$ is the distribution tangent to the $\rk$
orbits, such that for any $\t \in \rk$ and any  nonzero vector  $ v \in E_{i}$ the Lyapunov
exponent of $v$ is equal to $\chi_i (\t)$, i.e.
 $$
   \lim _{n \rightarrow \pm \infty }
   n^{-1} \log \| D(n\t) \, v \| = \chi_i (\t),
 $$
where $\| \cdot \|$ is any continuous norm on $TM$.
Any point $x\in \Re$ is called a {\em regular point}.

Furthermore, for any $\e > 0$ there exist positive measurable functions $C_\e (x)$
and $K_\e (x)$ such that for all $x \in \Re$, $v \in E_i(x)$, $\t \in \rk$, and $i=1, \dots, l$,
 \begin{enumerate}

\item $C^{-1}_\e (x) e^{\chi_i(\t)-\frac{1}{2} \e\|\t\|} \| v\|\le
\|D \t \,v \|  \le  C_\e (x) e^{\chi_i(\t)+\frac{1}{2} \e\|\t\|} \| v\|$;

\item Angles $\angle (E_{i}(x),T\O) \ge K_\e (x)$ and
         $\angle (E_{i}(x),E_{j}(x)) \ge K_\e (x)$, $i \not= j$;

\item $C_\e (\t x) \le C_\e (x) e^{\e \| \t \|}$
 and $K_\e (\t x) \ge K_\e (x) e^{-\e \| \t \|}$;

\end{enumerate}

\end{proposition}

The stable and unstable distributions of an element $\a(\t)$ are defined as
the sums of the Lyapunov distributions corresponding to
the negative and the positive Lyapunov exponents for $\at$ respectively:
$$
E^- _{\at} = \bigoplus_ {\chi_i (\t) <0} E_{i}, \qquad
E^+ _{\at} = \bigoplus_ {\chi_i (\t) >0} E_{i}.
$$

\subsection{Actions with Lyapunov exponents in general position.} \label{max rank}
Let $\a$ be an $\rk$ action as in the Main Theorem.
 Since  $(k+1)$ nontrivial Lyapunov exponents of $\a$ with
respect to $\mu$ are nonzero functionals and the Lyapunov hyperplanes are
in general position,  the  total number of Weyl chambers is equal to  $2^{k+1}-2$.
Each Weyl chamber corresponds to a different combination of signs for the
Lyapunov exponents. In fact, $2^{k+1}-2$ Weyl chambers correspond to all
possible combinations of signs except for all pluses and all minuses. The fact
that these two combinations are impossible can be seen as follows. First we
note that $\mu$ is non-atomic since it is ergodic for $\a$ and the entropy for some
element is positive. Now assume that there is an element $\t\in\rk$ such that all
exponents for $\a(\t)$ are negative. Then every ergodic component for $\a(\t)$ is
an isolated contracting periodic  orbit \cite[Proposition 1.3]{K97} and hence  the
measure $\mu$ must be atomic. In particular, we obtain the following property which
will play an important role in our considerations. Let $\chi_i$, $i=1,\dots,k+1$, be the
Lyapunov exponents of the action $\a$ and let $E_i$, $i=1,\dots,k+1$, be the
corresponding Lyapunov distributions.
\medskip

$(\mathcal C)$ {\em For every $i\in\{1,\dots,k+1\}$ there exists a Weyl chamber
$\mathcal C_i$ such that for every  $\t\in\rk\cap \mathcal C_i$ the signs of the
Lyapunov exponents are}
$$\chi_i(\t)<0\,\,  \text{and}\,\, \chi_j(\t)>0\,\,  \text{for all}\,\,  j\neq i.$$
In other words, property $(\mathcal C)$ implies that each Lyapunov distribution
$E_i$ is the full stable distribution for any $\t\in \mathcal C_i$.

Recall that stable distributions are always H\"older continuous (see, for example,
\cite{BPbook}).
Therefore, property $(\mathcal C)$ implies, in particular, that all Lyapunov
distributions for such actions inherit the H\"older continuity of stable distributions.
More generally, we have the following.

 \begin{proposition} \label{Holderness}
Let $\a$ be a $\Ce,\, \theta >0$ action  $\R^k$ as in Proposition \ref{MET}.
Suppose that a Lyapunov distributions $E$ is the intersection of the stable distributions
of some elements of the action. Then $E$ is H\"older continuous on any Pesin set
\begin {equation} \label{Pesinset}
\Re_{\e}^l =\{ x \in \Re  : \; C_\e (x) \le l, K_\e (x) \ge l^{-1} \}
\end{equation}
with H\"older constant which depends on $l$ and H\"older exponent $\delta>0$ which
depends on the action $\a$ only.

\end{proposition}

\subsection{Invariant manifolds}\label{sectionPesin}

We will use the standard material on invariant  manifolds corresponding to the
negative and positive Lyapunov exponents (stable and unstable manifolds) for $\Ce$
measure preserving diffeomorphisms of compact manifolds, see for example
\cite[Chapter 4]{BP}.

We will denote by $\wl_\at^-(x)$ and $\w_\at^-(x)$ correspondingly the local and
global stable manifolds for the diffeomorphism $\at$ at a regular point $x$.
Those manifolds  are tangent to the stable distribution $E_\at^-$. The global manifold is an
immersed Euclidean space and is defined uniquely. Any local manifold is a $\Ce$
embedded open disc in a Euclidean space. Its germ at $x$ is uniquely defined and
for any two choices their intersection is an open neighborhood of the point $x$ in each
of them. On any Pesin set $\Rel$ the local stable manifolds can be chosen of a
uniform size and changing continuously in the $\Ce$ topology.
The local and global unstable manifolds $\wl_\at^+(x)$ and $\w_\at^+(x)$
are defined as the stable manifolds for the inverse map $\a(-\t)$ and thus have
similar properties.

It is customary to use words ``distributions'' and ``foliations'' in this setting although
in fact the objects  we are dealing  with  are correspondingly measurable  families
of tangent spaces defined almost everywhere with respect to an invariant measure,
and  measurable families of smooth manifolds which fill a set of full measure.

Let $\a$ be an $\rk$ action as in the Main Theorem. Then property $(\mathcal C)$ 
gives that each Lyapunov distribution $E$ coincides with the full stable
distribution of some element of the action. Therefore, we have the corresponding
local and global manifolds $W (x)$ and $\w  (x)$ tangent to $E$. More generally,
these local and global manifolds are defined for any Lyapunov distribution $E$
as in Proposition \ref{Holderness}. We will refer to these manifolds as local and
global leaves of the Lyapunov foliation $\w$.




\section{Proof of Main Theorem}

As we mentioned before, part (1) of the Main Theorem follows immediately 
from part (2) by passing to the suspension. In this section we deduce part (2) 
from the technical Theorem \ref{tech}. First we show that the existence of an 
element with positive entropy implies that all regular elements have
positive entropy and that the conditional measures on every Lyapunov
foliation are non-atomic almost everywhere. This is done in Section
\ref{non-atomic}. Applying Theorem \ref{tech} we obtain that for
every Lyapunov foliation $\w$ the conditional measures on $\w$ are
absolutely continuous. We conclude the proof of the Main Theorem in
Section \ref{conclusion} by showing, as in \cite{KaK3}, that this
implies absolute continuity of $\mu$ itself.

\subsection{Conditional measures on all Lyapunov foliations are non\-atomic} \label{non-atomic}
We recall that a diffeomorphism has positive entropy with respect to
an ergodic invariant measure $\mu$ if and only if the conditional
measures of $\mu$ on its stable (unstable) foliation are non-atomic
a.e. This follows for example from \cite{LY2}. Thus if the entropy
$h_\mu (\t)$ is positive for some element $\t \in \rk$ then the
conditional measures of $\mu$ on $\w ^+_\at$ are non-atomic. Then
there exists an element $\s$ in a Weyl chamber $\mathcal C _i$ such
that the one-dimentional distribution $E_i=E^-_\as$ is not contained
in $E ^+_\at$ and thus $E ^+_\at \subset E^+_\as=\bigoplus_{j \ne i}
E_j$. Hence the conditional measures on $\w^+_\as$ are also
non-atomic. This gives  $h_\mu (\s)>0$ which implies that the
conditional measures on $\w_i=\w^-_\as$ must be also non-atomic. Now
for any $j \ne i$ consider the codimension one distribution $E'_j =
\bigoplus_{k \ne j} E_k =E ^+_{\a(\t_j)}$ for any $\t _j$ in the
Weyl chamber $\mathcal C _j$. Since $E_i \subset E'_j$ we see that
the conditional measures on the corresponding foliation $\w'_j$ are
non-atomic. Hence $h_\mu (\t_j)>0$ and the conditional measures on
$\w_j = \w^-_{\a(\t_j)}$ are non-atomic too. We conclude that the
conditional measures on every Lyapunov foliation $\w_i$, $i=1,
\dots, k+1$, are non-atomic. This implies, in particular, that the
entropy is positive for any non-zero element of the action.

\subsection{The absolute continuity of $\mu$} \label{conclusion}
The remaining argument is  similar to that in   \cite{KaK3}. In
order to prove that $\mu$ is an absolutely continuous measure we
shall use the following theorem that is essentially the flow
analogue of what is done in Section 5 of \cite{L} (see particularly
\cite[Theorem (5.5)]{L} and also \cite[Corollary H]{LY2}).
\begin{theorem}\label{srbimpliesleb}
Let $f:M\to M$ be a $C^{1+\alpha}$ diffeomorphism with invariant
measure $\mu$ and assume that $\h_{\mu}(f)$ is equal both to the sum
of the positive Lyapunov exponents and  to the absolute value of the
sum of the negative Lyapunov exponents. If the directions
corresponding to zero Lyapunov exponents integrates to a smooth
foliation and the conditional measures with respect to this central
foliation are absolutely continuous, then $\mu$ is absolutely
continuous with respect to Lebesgue measure.
\end{theorem}


To use Theorem \ref{srbimpliesleb} recall that there are $2^{k+1}-2$
Weyl chambers for $\a$ and any combination of positive and negative
signs for the Lyapunov exponents, except for all positive or all
negative, appears in one of the Weyl chambers. We use notations of
Section \ref{max rank} and consider an element $\t$ in the Weyl
chamber $\mathcal -C_i$. Then the Lyapunov exponents of $\t$ have
the following signs: $\chi_i (\t) >0$ and $\chi_j (\t)<0$ for all
$j\neq i$. Since the conditional measures on $\w ^+_{\at}$ are
absolutely continuous by Lemma~\ref{conditional abs cont}, we obtain
that
$$\h_{\mu}(\a(\t))=\chi_i(\t)$$
for any $\t$ in $- \mathcal C_i$.  By the Ruelle entropy inequality
$\h_{\mu}(\a(\t)) \le -\sum_{j\neq i}\chi_j(\t)$ and hence
$$\sum_{j=1}^{k+1}\chi_j(\t)\le0.$$
If $\sum_{j=1}^{k+1}\chi_j(\t)=0$ then Theorem~\ref{srbimpliesleb} applies  and the
proof is finished.

Thus we have to consider the case when
$\sum_{j=1}^{k+1}\chi_j(\t)<0$ for all $\t$ in all Weyl chambers
$-\mathcal C_i,\,\,i=1,\dots k+1$. This implies that
$\bigcup_{i=1}^{k+1}\mathcal C_i$ lies in the positive half space of
the linear functional $\sum_{j=1}^{k+1}\chi_j$. But this is
impossible since there exist elements $\t_i\in \mathcal
C_i,\,\,i=1,\dots k+1$  such that $\sum_{i=1}^{k+1}\t_i=0$. \qed


\section{The technical Theorem} \label{section tech}

In the notations of Proposition \ref{MET}, an ergodic invariant measure $\mu$
for a smooth locally free $\rk$ action $\a$ is called {\em hyperbolic} if all
nontrivial Lyapunov exponents  $\chi_i$, $i=1, \dots, l$, are nonzero
linear functionals on $\rk$.

\begin{theorem}  \label{tech}
Let $\mu$ be a hyperbolic ergodic invariant measure for a locally free
$\Ce,\, \theta >0$, action $\a$ of $\R^k,\,\,k\ge 2$, on a compact smooth
manifold $M$. Suppose that a Lyapunov exponent $\chi$ is simple and 
there are no other exponents proportional to $\chi$. Let $E$ be the
one-dimensional Lyapunov distribution corresponding to the exponent $\chi$.

Then $E$ is tangent $\mu$-a.e. to a Lyapunov foliation $\w$ and the conditional
measures of $\mu$ on $\w$ are either atomic a.e. or absolutely continuous a.e.
\end{theorem}

The assumptions on the Lyapunov exponents in Theorem~\ref{tech}
are considerably more general  than in the Main Theorem. In particular they may be satisfied 
for all exponents of a hyperbolic measure  for  an action on any rank greater than one   on a manifold of arbitrary large dimension. As an example one can take  restriction of an action satisfying the assumption of  part (1) of the Main Theorem to any lattice  $L\subset \Zk$
of rank at least  two which  has trivial intersection with all Lyapunov hyperplanes. 
For this reason Theorem~\ref{tech} has applications beyond the maximal rank case  considered in the Main Theorem. Those application with be discussed in a subsequent paper. 

On the other hand, positivity of entropy  for  some or even all non-zero elements is not sufficient to exclude atomic  measures on some of the Lyapunov foliations. Thus application to more general actions  may include stronger assumptions on ergodic properties of the measure.

\subsection{Outline of the proof of Theorem \ref{tech}} \label{outline tech}
We note that the Lyapunov distribution $E$ may not coincide with the full
stable distribution of any element of $\a$. First we will show that $E$
is an intersection of some stable distributions of $\a$.

An element $\t \in \rk$ is called {\em generic singular} if it belongs to exactly
one  Lyapunov hyperplane. We consider a generic singular element $\t$ in $L$,
i.e. $\chi$ is the only non-trivial Lyapunov exponent that vanishes on $\t$.
Thus
$$TM = T\O \oplus E_\at^- \oplus E \oplus E_\at^+$$
We can take a regular element $\s$ close to $\t$ for which $\chi (\s) > 0$ and
all other non-trivial  exponents have the same signs as for $\t$.
Thus $E_\as^-=E_\at^-$ and $E_\as^+=E_\at^+\oplus E$.
Similarly, we can take a regular element $\s'$ close to $\t$ on the other side
of $L$ for which $\chi (\s') < 0$ and $E_{\a(\s')}^+=E_\at^+$ and
$E_{\a(\s')}^-=E_\at^-\oplus E$. Therefore,
$$E = E_\as^+ \cap E_{\a(\s')}^- = E_{\a(-\s)}^- \cap E_{\a(\s')}^-.$$

We conclude that the Lyapunov distribution $E$ is an intersection of stable
distributions and, as in Proposition \ref{Holderness}, is H\"older continuous
on Pesin sets. As in Section \ref{sectionPesin}, $E$ is tangent $\mu$-a.e.
to the corresponding Lyapunov foliation $\w=\w_{\a(-\s)}^- \cap \w_{\a(\s')}^-$.

We denote by $\mw_x$ the system of conditional measures of $\mu$ on $\w$.
By ergodicity of $\mu$ these conditional measures
are either non-atomic of have atoms for $\mu$-a.e. $x$. Since $\w$ is an invariant
foliation contracted by some elements of the action, it is easy to see that in the latter
case the conditional measures are atomic with a single atom for $\mu$-a.e. $x$
(see, for example, \cite[Proposition 4.1]{KS3}).
The main part of the proof is to show that if the conditional measures $\mw_x$ are
non-atomic for $\mu$-a.e. $x$, then they are absolutely continuous $\mu$-a.e.

To prove absolute continuity of the conditional measures on $\w$ we show in
Section \ref{sectionconclusion} that they are Haar with respect to the invariant
family of smooth affine parameters on the leaves of $\w$.
As in \cite{KaK3}, this uses affine maps of the leaves which
preserve the conditional measures up to a scalar multiple. Such affine maps are
obtained in \cite{KaK3} as certain limits of actions along $\w$ by some elements
of the action. It is essential that derivatives of these elements along $\w$
are uniformly bounded. In \cite{KaK3} it was possible to choose such elements
within the Lyapunov hyperplane $L$. We note that in general  Multiplicative Ergodic
Theorem  only guarantees that the elements in $L$ expand or contact  leaves of
$\w$ at a subexponential rate.

The main part of the proof is to produce a sequence of elements of the
action with uniformly bounded derivatives along $\w$ and with enough recurrence. In
Section \ref{Lyapunov metric} we define a special Lyapunov metric on distribution $E$
and show that it is H\"older continuous on Pesin sets. Then in Section \ref{measurable
time change} we construct a {\em measurable} time change for which the expansion
or contraction in $E$ with respect to this Lyapunov metric is given exactly by the
Lyapunov exponent $\chi$. This gives sufficient control for the derivatives along $\w$.

To produce enough recurrence we study  properties of this measurable time
change in Section \ref{time change properties}. We prove that it is differentiable along
regular orbits and H\"older continuous when restricted to  any Pesin set \eqref{Pesinset}.
This allows us to
show in Section \ref{new action properties} that the time change has some structure
similar to that of the original action. First, it preserves a measure equivalent to $\mu$.
Second, it preserves certain  ``foliations'' whose restrictions to Pesin sets are H\"older
graphs over corresponding foliations of the original action $\a$. More precisely, the
leaves for $\a$ are tilted along the orbits  to produce invariant sets for the time change
action $\beta$ and the tilt is a H\"older function when restricted to a set of large measure
(the intersection with such a set has large conditional measure for a typical leaf).
Of course, the H\"older constants (but not the exponents) deteriorate when one increases
the Pesin set but in the end one gets a measurable function defined almost everywhere.

Using these properties, we show in Section \ref{ergodicitysingular} that for a typical
element in the Lyapunov hyperplane $L$ the time change acts sufficiently transitively
along the leaves of $\w$. For this we use the ``$\pi$-partition trick'' first introduced in
\cite{KS3} for the study of invariant measure of actions by automorphisms of a torus and
adapted to the general nonuniform situation in \cite{KaK3}. We use this argument for the
time change action $\beta$ and the main technical difficulty is in showing that the weird
``foliations'' described above can still be used in essentially the same way as for smooth
actions.


\section{Lyapunov metric} \label{Lyapunov metric}

In this section we use notations of Theorem \ref{tech}.
We define a Lyapunov metric on the Lyapunov distribution $E$
and establish its properties.

We fix a smooth Riemannian metric $<\cdot,\cdot>$ on $M$. Given $\e>0$
and a regular point $x \in M$ we define {\em the standard $\e$--Lyapunov scalar
product (or metric)} $<\cdot,\cdot>_{x,\e}$ as follows. For any $u,v\in E(x)$ we define
\begin{equation}\label{eqLyapunovmetric}
<u,v>_{x,\e}=\int_{\rk}<(D \s)  u, (D \s)  v>
\exp(-2\chi(\s)-2\e\|\s\|) \, d\s
\end{equation}
We observe using (1) of Proposition \ref{MET} that for any $\e >0$
the integral above converges exponentially for any regular point $x$.

We will usually omit the word ``standard'' and will call this scalar product
{\em $\e$--Lyapunov metric} or, if  $\e$ has been fixed and no confusion
may appear, simply {\it Lyapunov metric}. The norm generated by this scalar
product will be called  the (standard $\e$--) {\em  Lyapunov norm} and denoted by
$\|\cdot\|_{x,\e}$ or $\|\cdot\|_{\e}$.

\begin{remark}
The definition above gives a measurable scalar product on any Lyapunov
distribution $E$ of an arbitrary nonuniformly hyperbolic $\rk$ action (and
similarly for a $\Zk$ action), without any assumption on Lyapunov exponents,
such as multiplicity, or on geometry of Lyapunov hyperplanes. One can also
define the Lyapunov scalar product on the whole tangent space $T_x M$ by
declaring  the Lyapunov distributions to be pairwise orthogonal and orthogonal
to the distribution $T\O$ tangent to the orbits of the $\rk$ action (on $T\O$ one
can take a canonical  Euclidean metric given by the action). Proposition
\ref{Lyapunov metric derivative} as well as estimates \eqref{Lyapunov>metric}
and \eqref{Lyapunov<metric} hold for such general case. Also, continuity of the
Lyapunov scalar product on sets of large measure follows simply from
measurability by Luzin's theorem. However, H\"older continuity on Pesin sets
for the Lyapunov scalar product on a given Lyapunov distribution $E$ requires
similar H\"older continuity of $E$. The latter is not necessarily true for an
arbitrary Lyapunov distribution.
\end{remark}

We denote by $D_x^{E}$ the restriction of the derivative to the Lyapunov
distribution $E$. The main motivation for introducing  Lyapunov metric is
the following estimate for the norm of this restriction with respect to the
Lyapunov norm.

\begin{proposition} \label{Lyapunov metric derivative}
For any regular point $x$ and any $\t \in \rk$
\begin{equation}
\exp(\chi (\t)-\e\|\t\|) \le
\|D_x^{E} \, \, \t \|_{\e}  \le   \exp(\chi (\t)+\e\|\t\|).
\label{Lmetric}
\end{equation}

\end{proposition}

\proof Fix a nonzero $u\in E(x)$. By the definition of the standard
$\e$--Lyapunov norm  we obtain
$$
\|(D_x \t) u\|^2_{\t x,\e}  =
\int_{\rk} \|(D_{\t x} \s)  (D_x \t) u \|^2  \exp(-2\chi(\s)-2\e\|\s\|) \, d\s =
$$
$$
\int_{\rk} \|(D_{x} (\s + \t)) u \|^2  \exp(-2\chi(\s)-2\e\|\s\|) \, d\s =
$$
\begin{equation}
\int_{\rk} \|(D_{x} \s') u \|^2  \exp(-2\chi(\s'-\t)-2\e\|\s'-\t\|) \, d\s'
\label{norm^2}
\end{equation}

We note that the exponent can be estimated above and below as follows
$$
     -2\chi(\s'-\t)-2\e\|\s'-\t\| \le  (-2\chi(\s')-2\e\|\s'\|) +2(\chi(\t) + \e\|\t\|),
$$
$$
     -2\chi(\s'-\t)-2\e\|\s'-\t\| \ge   (-2\chi(\s')-2\e\|\s'\|) +2(\chi(\t) - \e\|\t\|).
$$
These inequalities together with the definition
$$
\|u\|_{x,\e}^2 = \int_{\rk} \|(D_{x} \s') u \|^2  \exp(-2\chi(\s')-2\e\|\s'\|) \, d\s'
$$
give the following estimate
$$
e^{2(\chi(\t)-\e\|\t\|)}\|u\|_{x,\e}^2 \le
\|(D_x \t) u\|_{\t x,\e}^2  \le   e^{2(\chi(\t)+\e\|\t\|)}\|u\|_{x,\e}^2
$$
which concludes the proof of  the proposition. \qed
\medskip

Now we establish some important properties of the Lyapunov metric.
First we note that the original smooth metric gives a uniform below
estimate for the Lyapunov metric, i.e. there exists positive constant $C$
such that for all regular $x \in M$ and all $u \in E$
\begin{equation}\label{Lyapunov>metric}
 \|u\|_{x,\e}\ge C \|u\|
\end{equation}

The next proposition establishes the opposite inequality as well as
continuity of the Lyapunov metric  on a given Pesin set. We note
that similarly to the proof of  Lemma~\ref{F} below one can
show that $\e$-Lyapunov metric is actually smooth along the orbits.

\begin{proposition} \label{Lyapunov metric continuity}
The $\e$-Lyapunov metric is continuous along any regular orbit and on any
Pesin set $\Rel$. Furthermore, there exists $C(l,\varepsilon) >0$ such that
for all $x \in \Rel$ and all $u \in E$
\begin{equation}\label{Lyapunov<metric}
\|u\|_{x,\e} \le C(l,\varepsilon) \|u\|.
\end{equation}
\end{proposition}

\proof
Let us fix $u \in E (x)$ with $\|u\|=1$. The integrand in
 equation \eqref{eqLyapunovmetric}
$$
f(x,\s) = <(D_x \s)  u, (D_x \s)  u>  \exp(-2\chi(\s)-2\e\|\s\|)
$$
is continuous with respect to $x$ on $\Rel$ by Proposition \ref{Holderness}. Also, by (1)
of Proposition \ref{MET} we have $|f(x,\s)| \le C_{\e}(x) \exp(-\e \|\s\|)$ and hence  for $x \in \Rel$
$$\int_{\rk} f(x,\s) \, d\s \le \int_{\rk} C_{\e}(x) \exp(-\e \|\s\|) \, d\s \le l \int_{\rk}  \exp(-\e \|\s\||) \, d\s.$$

This implies the estimate \eqref{Lyapunov<metric} and the continuity of the metric
on the Pesin set $\Rel$. The continuity along orbits follows since for any regular point
$x$ and any bounded set $B \subset \rk$ there is $l$ such that $Bx \subset \Rel $.
\qed\medskip

Next we obtain H\"older continuity of the Lyapunov metric which will be crucial
for deducing   properties of the time change in Section \ref{TC}.

\begin{proposition} \label{Lyapunov metric Holder}
There exists $\gamma >0$ which depends only on $\e$ and on the action
and $K(l,\epsilon) > 0$ which, in addition, depends on the Pesin set $\Rel$
such that the $\e$-Lyapunov metric is H\"older continuous on $\Rel$ with
exponent $\gamma$ and constant $K(l, \epsilon)$.
\end{proposition}

\begin{remark} We note the dependence of the constant in \eqref{Lyapunov<metric}
and H\"older constants in Propositions \ref{Holderness} and  \ref{Lyapunov metric Holder},
as well as below in  Propositions  \ref{time change Holder} and \ref{new stable} on the
Pesin set $\Rel$. For a fixed $\e$ these constants depend only on $l$ and can be
estimated by $Cl^p$ for some power $p$. This holds in Proposition \ref{Holderness}
and can be observed in the proofs of the other propositions.
By (3) of Proposition \ref{MET} for any $\t \in \rk$ we have $\at (\Rel) \subset \Re_\e ^{l'}$
with $l'=\exp (\e \| \t \|) l$. Therefore, we can say that these constants may grow in $\t$
with a slow exponential rate, more precisely, by a factor at most $\exp (p \e \| \t \|)$.
\end{remark}

\proof
We consider two
nearby points $x$ and $y$ in a Pesin set $\Rel$. By Proposition \ref{Holderness} we
can take vectors $u \in E(x)$ and $v \in E(y)$ with $\|u\|=\|v\|=1$ for which the distance
in $TM$ can be estimated as $\dist (u,v) \le K_1 \rho ^{\d}$, where $\rho = \dist (x,y)$.
Since $E$ is one-dimensional it suffices to show that $\|u\|_{x,\e}$ and $\|v\|_{y,\e}$
are H\"older close in $\rho$. We will now estimate
$$
\|u\|_{x,\e}^2 - \|v\|_{y,\e}^2=\int_{\rk}( \|(D \s) u\|^2 - \|(D \s)  v \|^2 )
\exp(-2\chi (\s)-2\e\|\s\|) \, d\s
$$
Using spherical coordinates $\s = \s (r, \theta)$ where $r=\|\s\|$ and denoting
$$
\psi (\s) = ( \|(D \s) u\|^2 - \|(D \s)  v \|^2 ) \exp(-2\chi (\s)-2\e r)
$$
we can write
$$
\|u\|_{x,\e}^2 - \|v\|_{y,\e}^2=\int_0^\infty f(r)\,dr ,
\;\; \text{where} \;\;
f(r)=r^{k-1} \int_{S^{k-1}} \psi (\s) \, d\theta .
$$
We will estimate the difference $\|(D \s) u\|^2 - \|(D \s)  v \|^2$
inside a large ball using closeness of $u$ and $v$ and outside it by estimating
each of the two terms. Since the action $\a$ is smooth we observe that
$$
|  \|(D \s) u\|^2 - \|(D \s)  v \|^2 | \le
\dist ((D \s) u , (D \s)  v) \cdot (\|(D \s) u\| + \|(D \s)  v \|)
$$
$$
 \le C_1 \exp (L\| \s \|) \cdot \dist (u , v) \le  C_1 \exp (Lr) \; K_1 \rho ^{\d}
$$
for some $L>0$. Hence we obtain $| \psi (\s)| \le K_1 C_2 \exp (M'r) \; \rho ^{\d}$,
where $M'=L+2(\| \chi \| + \e)$. Then for sufficiently large $a$ we have
$$
\int_0^a |f(r)| \, dr
\le \int_0^a  r^{k-1} K_1 C_3 \exp ((M')r) \, \rho ^{\d} dr  \le
K_1 C_4 \exp (Ma) \, \rho ^{\d},
$$
where for simplicity we absorbed the polynomial factor appearing in the estimates
into the exponent. Then
\begin{equation}\label{spherical}
\text{for }  \quad a=\frac{\d}{2M} \log \frac1{\rho} \quad \text{ we have }
\int_0^a  |f(r)| \, dr  \le K_2  \, \rho ^{\d/2}. \end{equation}
Now we consider $\int_a^\infty |f(r)| \, dr$. Since $x,y \in \Rel$, using (1) of
Proposition \ref{MET} we obtain $| \psi (\s) | \le 2 l^2 \exp (-\e r)$.
Hence
$$
\int_a^\infty |f(r|) \, dr
\le \int_a^\infty  r^{k-1} C_5 l^2 \exp (-\e r) \, dr  \le
K_3 \exp (-\e a/2)
$$
where we again absorbed the polynomial factor into the exponent.
For $a$ defined in \eqref{spherical} this gives us
$$
\int_a^\infty |f(r)| \, dr \le K_3 \exp (-\e a/2) \le K_3 \rho^{\gamma} ,
$$
where $\gamma=\frac{\e \d}{4M}$. Combining this with the estimate \eqref{spherical}
for $\int_0^a |f(r)| \, dr$ we obtain
$$
|\|u\|_{x,\e}^2 - \|v\|_{y,\e}^2| \le \int_0^\infty |f(r)| \,dr \le K_4 \rho^{\gamma}.
$$
According to \eqref{Lyapunov>metric}, the Lyapunov norm is bounded below
by the usual norm,  so that
$$
|\|u\|_{x,\e} - \|v\|_{y,\e}| \le |\|u\|_{x,\e}^2 - \|v\|_{y,\e}^2|/(\|u\|_{x,\e} + \|v\|_{y,\e}) \le
$$
$$
\le |\|u\|_{x,\e}^2 - \|v\|_{y,\e}^2|/2K \le K_5 \rho^{\gamma}
$$
which completes the desired H\"older estimate. \qed


\section{Measurable  time change and its properties}\label{TC}


\subsection{Construction of a measurable time change} \label{measurable time change}

In this section we use notations of Theorem \ref{tech}.
We fix small $\e >0$ and consider the Lyapunov metric $\|.\|_\e$ on the Lyapunov
distribution $E$. We first study the behavior of the derivative restricted to $E$ along the
$\rk$-orbits. For a regular point $x$ we consider the function $f_x : \rk \to \R$ given by
\begin{equation}
  f_x (\t) = \log \| \de_x \t \|_{\e}
\label{Fdef}
\end{equation}
According to Proposition \ref{Lyapunov metric derivative} the function $f$ satisfies
inequalities
\begin{equation}
\chi (\t)-\e\|\t\| \le
f_x ( \t )  \le  \chi (\t)+\e\|\t\|.
\label{Lmetricf}
\end{equation}
Also, since $E$ is one-dimensional, $f$ satisfies the cocycle identity
\begin{equation}
 f_x (\t + \s) = f_x (\t) + f_{\t x} (\s).
\label{Fcocycle}
\end{equation}

We will now establish smoothness of the function $f_x$ in $\t$.

\begin{lemma}  \label{F}
For any regular point $x$ the function $f_x (\t)$ is $C^1$.
More precisely, for any $\t,\ev \in \rk$ we have
$$
(D_\t f_x )\, \ev = \chi (\ev) + \e \, \psi _{\t x}(\ev),
$$
where $|\psi _{\t x}(\ev)|\le \frac12 \| \ev \|$ and $\psi _{\t x}(\ev)$ is continuous
in $\t$ and $\ev$.

\end{lemma}

\proof
Fix a regular point $x$  and
 consider the function
$$
F (\t)= \exp (f_x(\t)) =\| \de_x \t \|_{\e}.
$$
Fix a vector $u \in E(x)$ with $\| u\|_{x,\e}=1$. Since $E(x)$ is one-dimensional,
$$
F (\t) = \|(D_x \t) u\|_{\t x,\e}.
$$
Using the definition of the Lyapunov metric we obtain as in \eqref{norm^2} that
$$
F^2 (\t) = \int_{\rk} \|(D_{x} \s) u \|^2  \exp(-2\chi(\s)+2\chi(\t)-2\e\|\s-\t\|) \, d\s
$$
Differentiating at $\t$ we obtain
$$
(D_\t \,F^2)\, \ev = \int_{\rk} \|(D_{x} \s) u \|^2
\exp(-2\chi(\s)+2\chi(\t)-2\e\|\s-\t\|) \times
$$
$$
 \left(2\chi (\ev) + \e \frac{<\s-\t,\ev>}{\|\s-\t \|} \right) \, d\s =
$$
$$
\int_{\rk} \|(D_{x} (\s + \t)) u \|^2  \exp(-2\chi(\s)-2\e\|\s\|)
\cdot \left(2\chi (\ev) + \e \frac{<\s,\ev>}{\|\s \|} \right) \, d\s =
$$
$$
= 2\chi (\ev) F^2(\t) + \e \tilde \psi (\t,\ev),
$$
where
$$
\tilde \psi (\t,\ev) = \int_{\rk} \|(D_{x} (\s + \t)) u \|^2  \exp(-2\chi(\s)-2\e\|\s\|)
 \frac{<\s,\ev>}{\|\s \|}  \, d\s
$$
Then for the function $f_x$ we obtain
$$
(D_\t f_x)\, \ev = \frac12 D_\t ( \log{F^2})\, \ev = \frac{(D_\t F^2)\, \ev}{2 F^2(\t)}=
\chi (\ev)  + \e \, \psi _x (\t,\ev),
$$
where
$\psi (\t,\ev) = \tilde \psi (\t,\ev) /2F^2(\t)$.
We observe that $\tilde \psi (\t,\ev)$ is continuous in $\t$ and $\e$ and hence
so is $\psi _x (\t,\ev)$. We conclude that $f_x(\t)$ is $C^1$.
Since $ |<\s,\ev> \cdot \|\s \| ^{-1}| \le \| \ev \| $ we obtain
$|\tilde \psi (\t,\ev)| \le F^2(\t ) \| \ev \|$ and hence
$$|\psi _x (\t,\ev)| \le \| \ev \|/ 2.$$
We also note that $\psi _x (\t,\ev)=\psi _{\t x}(\0,\ev)$, which follows for example
from the cocycle relation \eqref{Fcocycle}.  Denoting $\psi _{x}(\ev)= \psi (\0,\ev)$
we obtain the desired formula for $(D_\t f_x)\, \ev$ with function $\psi _{x}(\ev)$
which is continuous in $\ev \in \rk$ and satisfies $|\psi _{\t x}(\ev)|\le \frac12 \| \ev \|$.
\qed
\medskip

Now we proceed to constructing the time change. We fix a vector $\nor$ in $\rk$
normal to $L$ with $\chi(\nor)=1$. We will assume that $\e$ and $\e \| \nor \|$ are
both small, in particular $\e \| \nor \|<1/2$.

\begin{proposition}  \label{adjust}
For $\mu$-a.e. $x\in M$ and any $\t \in \rk$ there exists a unique
real number $g(x,\t)$ such that the function  $\g (x,\t)=\t+g(x,\t) \nor$
satisfies the equality  $$\| \de _x \a (\gxt) \|_\e = e^{\chi (\t)}.$$
The  function $g(x,\t)$ is measurable and is continuous in $x$ on Pesin sets  \eqref{Pesinset} and
along the orbits of $\a$. It satisfies the inequality  $|g(x,\t)| \le 2\e \| \t \|$.
\end{proposition}

In Section \ref{time change properties} we will show that $g(x,\t)$
is actually H\"older continuous in $x$ on Pesin sets and is  $C^1$ in $\t$.

\proof
Recall that by Proposition \ref{Lyapunov metric derivative} for any regular
point $z$ we have
\begin{equation}
\exp (\chi (\t) - \e\| \t \|) \le \|\de _z \at \|_\e \le \exp (\chi (\t) + \e\| \t \|);
\label{adjust1}
\end{equation}
thus  in particular
\begin{equation}
\exp (s - \e s \| \nor \|) \le \|\de _z \a(s \nor) \|_\e \le \exp (s + \e s \| \nor \|).
\label{adjust2}
\end{equation}
We fix a regular point $x$ and define
$$
\phi(s) =\log \|\de _{\at x} \a(s \nor) \|_\e =f_{\at x} (s \nor).
$$
Using Lemma \ref{F} we obtain
\begin{equation}
\phi '(s) = \frac{d}{ds} f_{\at x} (s \nor)=(D_{s\nor} f_{\at x}) \nor=
\label{adjust3}
\end{equation}
$$
= \chi (\nor)  \, + \e \, \psi _{(s\nor +\at) x} (\nor)) \ge  1- \e  \| \nor \| /2 >0
$$
This implies that $\phi: \R \to \R$ is a $C^1$  bijection.
Hence there exists a unique number $s_0$ such that
$\phi(s_0)= \chi (\t) - \log \|\de _x \at \|_\e$ and thus $g (x,\t)=s_0$
satisfies the equation in the lemma. We observe that \eqref{adjust1} implies
$$- \e\| \t \| \le \chi (\t) - \log \|\de _x \at \|^{-1}_\e \le \e\| \t \|.$$ Also, \eqref{adjust2}
implies $$s - \e s \| \nor \| \le \phi (s) \le s + \e s \| \nor \|.$$ Hence $|\phi (s)| \ge
\frac12 |s|$ and we conclude that $|g (x,\t)| =|s_0| \le 2\e\| \t \|$

The continuity of $g (x,\t)$ in $x$ on Pesin sets and along the orbits of $\a$
follows from the corresponding continuity of the Lyapunov norm. \qed

\begin{proposition}  \label{newaction}
The formula $\b (\t,x)= \a (\gxt)x$ defines an $\rk$ action $\b$ on $M$
which is a measurable time change of $\a$, i.e.
\begin{equation}
\b (\s +\t, x) = \b (\s, \b (\t ,x)), \qquad \text{or}
\label{b-action}
\end{equation}
\begin{equation}
\g (x, \s +\t )= \gxt + \g (\a (\gxt)x,\s).
\label{g-timechange}
\end{equation}
The action $\b$ is measurable and is  continuous on any Pesin set for $\a$.
\end{proposition}

\begin{remark} The time change is defined using a condition on the derivative
of the {\em original} action restricted to $E$. The new action is not necessarily
smooth and  typically does not  preserve the Lyapunov foliation $\w$. However,
it does preserve the sum of the distribution $E$ with the orbit distribution as well
as the corresponding {\em orbit-Lyapunov} foliation.
\end{remark}

\proof
We will verify \eqref{b-action}. This relies on the uniqueness part of the previous
proposition. If we denote
$$y =\b (\t,x)= \a (\gxt)x $$
we can rewrite the right side of \eqref{b-action} as

\begin{equation}\begin{aligned} \b (\s, \b (\t ,x)) =  \b (\s, y) = & \a (\g (y,\s))\,
y = \\
\a (\g (y,\s)) \circ \a (\gxt) \,x =
\a &(\g (y,\s) +\gxt ) \, x = \\
 \a (\s + \t + (g (y,\s) &+g(x,\t)) \nor) \, x. \label{b}\end{aligned}
\end{equation}

From this equation we see that the point $\b (\s, \b (\t ,x))$ belongs to
the $\{t \nor \}$-orbit of $\a (\s +\t) x$.  By definition,
the point $\b (\s +\t, x)$ also belongs to this orbit, moreover, it is the
unique point of the form $\a (\s +\t +g \nor) x$ for which
$$
\| \de _x \a (\s +\t +g \nor) \|_\e = e^{\chi (\s +\t)}.
$$
On the other hand, equation \eqref{b} and the definition of $\g (y,\s)$ and $\gxt$
imply that
$$\begin{aligned}
\| \de _x \a (\s + \t + (g (y,\s) &+g(x,\t)) \nor) \|_\e =\\
\| \de _y \a (\g (y,\s)) \|_\e  \cdot  \| \de _x \a& (\gxt)\|_\e = e^{\chi (\s)} \cdot e^{\chi (\t)}.\end{aligned}
$$
Thus we conclude that the points $\b (\s, \b (\t ,x))$
and $\b (\s +\t, x)$ coincide, i.e. \eqref{b-action}. In particular, we obtain
$$
g (x, \s +\t )= g(x,\t) + g (\b (\t,x),\s)
$$
which gives equation \eqref{g-timechange}.
\qed

\subsection{Properties of the time change} \label{time change properties}

\begin{proposition} \label{time change Holder}
The time change $\g (x,\t)$ is H\"older continuous in $x$ with H\"older
exponent $\gamma$ on any Pesin set $\Rel$. The H\"older constant depends
on the Pesin set and can be chosen uniform in $\t$ for any compact subset of $\rk$.
\end{proposition}

\proof
We fix $\t \in \rk$ and two nearby points $x$ and $y$ in a Pesin set $\Re_\e ^{l'}$.
We take $l=l' \exp (\e \| \t \|)$ and note that by Proposition \ref{MET} (3) we have
$\at (\Re_\e ^{l'}) \subset \Rel$. Hence the points $x$, $y$, $\a(\t)x$, and
$\a(\t)y$ are all in the Pesin set $\Rel$. To prove the H\"older continuity of $\g (x,\t)$
we need to show that $|g(x,\t)-g(y,\t)|$ are H\"older close with respect to the distance
between $x$ and $y$, i.e. can be estimated from above by a constant multiple
of a power of  $\rho=\dist(x,y)$.

First we show that
$\| \de _x \a (\t+g(x,\t) \nor) \|_\e$ and $\| \de _y \a (\t+g(x,\t) \nor) \|_\e$
are H\"older close in $\rho$. This can be seen as follows. Since the  action $\a$
is smooth, the points $\a (\t+g(x,\t) \nor,x)$ and $\a (\t+g(x,\t) \nor,y)$ as well
as the derivatives $D _x \a (\t+g(x,\t) \nor)$ and $D _y \a (\t+g(x,\t) \nor)$ are
H\"older close
in $\rho$ with constant depending only on the action and $\| \t \|$.
Also, by Proposition \ref{Holderness} the distribution $E$ is H\"older continuous in $\rho$
on the Pesin set $\Rel$. Finally, the Lyapunov metric is H\"older continuous in $\rho$
on $\Rel$ with the H\"older exponent $\gamma$ by Proposition \ref{Lyapunov metric
Holder} and its ratio to a smooth metric is uniformly bounded above and below on
$\Rel$ by \eqref{Lyapunov>metric} and \eqref{Lyapunov<metric}. We conclude that
\begin{equation}
| \| \de _x \a (\t+g(x,\t) \nor) \|_\e - \| \de _y \a (\t+g(x,\t) \nor) \|_\e | \le K_1 \rho^\gamma
\label{gH1}
\end{equation}
By the definition of $\g (x,\t)=\t+g(x,\t) \nor$ we have
\begin{equation}
\| \de _x \a (\t+g(x,\t) \nor) \|_\e = e^{\chi (\t)} =\| \de _y \a (\t+g(y,\t) \nor) \|_\e .
\label{gH2}
\end{equation}
Then \eqref{gH1} and the first equality in \eqref{gH2} imply that
$$| \| \de _y \a (\t+g(x,\t) \nor) \|_\e - e^{\chi (\t)}| \le K_1 \rho^\gamma.$$
We note that the points $\a (\t+g(x,\t)\nor,y)$ and $\a (\t+g(y,\t)\nor,y)$ are
on $\{t \nor \}$-orbit of point $\a (\t ,y)$, and that the value $g(y,\t)$ is determined
by the second equality in \eqref{gH2}. Therefore, the difference $g(y,\t)-g(x,\t)$
represents the time adjustment in $\nor$ direction required to bring the norm
$\| \de _y \a (\t+s \nor) \|_\e$
from being $K_1 \rho^\gamma$-close to $e^{\chi (\t)}$ to being exactly $e^{\chi (\t)}$.
Recall that by Lemma~\ref{F} the norm $\| \de _y \a (\t+s \nor) \|_\e $ varies smoothly
with $s$ (see equation \eqref{adjust3} in the proof of Proposition \ref{adjust}).
Thus we conclude that
$|g(y,\t)-g(x,\t)|\le K_2 \rho^\gamma$.
\qed

\begin{proposition}  \label{time change smooth}
The time change $\g (x,\t)=\t+g(x,\t) \nor$ is differentiable and $C^1$
close to identity in $\t$. More precisely, for a.e. $x$
$$
 \| \frac{\partial g} {\partial \t} (x,\t) \| \le \e .
$$

\end{proposition}

\proof
We fix a regular point $x$ and vectors $\t ,\ev$ in $\rk$ and
consider a function of two real variables
\begin{equation}
\Phi (s,g) = f_x (\t +s \ev + g \nor) - \chi (\t +s \ev).
\label{Phidef}
\end{equation}
We note that, by Proposition \ref{adjust}, $g(s)=g(x,\t +s \ev)$ is the unique
solution for the implicit function equation $\Phi (s,g(s)) = 0$.

Using Lemma \ref{F} we obtain
$$
\frac{\partial \Phi} {\partial g} (s,g) = (D_{\t +s \ev + g \nor} f_x )\, \nor
= \chi (\nor) \, + \e \, \psi _{(\t +s \ev + g \nor)x} (\nor),
$$
where $|\psi _{(\t +s \ev + g \nor)x}(\nor)|\le \frac12 \| \nor \|$ and
 is continuous in $s$ and $g$. Similarly,
$$
\frac{\partial \Phi} {\partial s} (s,g)  = \chi (\ev) \,
+ \e \, \psi_{(\t +s \ev + g \nor)x} (\ev) -\chi (\ev)=\e \, \psi_{(\t +s \ev + g \nor)x} (\ev).
$$
We conclude that $\Phi$ is a $C^1$ function of $(s,g)$. Moreover,
since  $\chi (\nor) =1$ and $\e$ is small we obtain
$$
\frac{\partial \Phi } {\partial g} (s,g)  = 1 + \e \, \psi  _{(\t +s \ev + g \nor)x} (\nor)
\ge 1- \e \| \nor \| /2 >0,
$$
Therefore, by implicit function theorem, $g(s)$ is differentiable and
$$
g'(s)= - \left(  \frac{\partial \Phi} {\partial s} (s,g(s)) \right) \left(
\frac{\partial \Phi } {\partial g} (s,g(s)) \right) ^{-1} =
\frac{-\e  \psi _{(\t +s \ev + g \nor)x} (\ev)}{1 + \e \, \psi _{(\t +s \ev + g \nor)x} (\nor)}.
$$
Moreover, since $|\psi _{(\t +s \ev + g \nor)x}(.)|\le \frac12 \| . \|$ we obtain
$$
|g'(s)| \le \frac{\e \| \ev \|}{2 - \e \, \| \nor \|} \le \e \| \ev \|
$$
provided that $\e \| \nor \|<1$. Since $g(s)=g(x,\t +s \ev)$ we have
$$
\left( \frac{\partial g} {\partial \t} (x,\t) \right) \ev = g'(s)
$$
and thus the partial derivatives of $g(x,\t)$ in the second variable
exist and are continuous in $\t$. We conclude that $g(x,\t)$ is $C^1$
in $\t$ with
$
 \| \frac{\partial g} {\partial \t} (x,\t) \| \le \e $. \qed

\subsection{Properties of the action $\b$.} \label{new action properties}
We note that the new action $\b$ is not smooth. Hence the notions and
results of nonuniformly hyperbolic theory do not apply to $\b$ formally.
In particular, such objects as derivatives, Lyapunov distributions,
Lyapunov exponents, and Lyapunov hyperplanes will always refer
to the ones of the original action $\a$. However, the new action $\b$
inherits some structures of $\a$ such as invariant measure and invariant
``foliations" which are close to those of $\a$. This is described
in the following two statements. We will use these structures in the
next section to obtain some important transitivity properties of $\b$.

\begin{proposition}  \label{new measure}
The determinant of the time change $\g (x,\t)$
$$\Delta (x) =\det \left( \frac{\partial \g}{\partial \t}  (x,\0) \right) $$
is a measurable function which is $L^\infty$ close to constant $1$ on $M$.
Therefore, the new action $\b$ preserves an invariant measure $\nu$
which is absolutely continuous with respect to $\mu$ (and equivalent to $\mu$)
with density $\frac{d\nu}{d\mu} = \Delta (x)^{-1}$.
\end{proposition}

\proof
The $L^\infty$ estimate for the determinant follows immediately from the
fact that by Proposition \ref{time change smooth} for a.e. $x$
$$\| \frac{\partial \g (x,\t)}{\partial \t} - \Id \| \le \e.$$
Then the existence of the invariant measure $\nu$ for $\b$ follows from \cite{Katok78}.
\qed
\medskip

We will denote by $\on$ the orbit foliation of the one-parameter subgroup $\{t \nor \}$.

\begin{proposition}  \label{new stable}
For any element $\s \in \rk$ there exists stable ``foliation" $\wb ^-_\bs$ which is
contracted by $\bs$ and invariant under the new action $\b$. It consists of ``leaves"
$\wb ^-_\bs (x)$ defined for almost every $x$. The ``leaf" $\wb ^-_\bs(x)$ is a measurable
subset of the leaf $(\on \oplus \w^-_\as)  (x)$ of the form
$$
\wb ^-_\bs (x) =\{ \a (\varphi_x (y) \nor) y : y \in \w^-_\as(x) \},
$$
where $\varphi_x : \w^-_\as(x) \to \R$ is an almost everywhere defined
measurable function. For $x$ in a Pesin set the $\varphi_x$ is H\"older
continuous on the intersection of this Pesin set with any ball of fixed radius
in $\w^-_\as(x)$with H\"older exponent $\gamma$ and constant which
depends on the Pesin set and radius.
\end{proposition}

\proof  We will give an explicit formula for the function $\varphi_x$ in terms of the
time change so that its graph is contracted by $\bs$. The calculation is similar to
finding stable manifolds for a time change of a flow. The H\"older continuity of
$\varphi_x$ will follow from the formula and the H\"older continuity of the time change.
Since $\w^-_\as$ is invariant under $\a$, we note that  $\on \oplus \w^-_\as$ is invariant
under $\b$ by the construction of the time change. Since $\wb ^-_\bs (x)$ is clearly
characterized within $\on \oplus \w^-_\as (x)$ by the contraction property and since
$\bt$ is continuous on Pesin sets, the usual argument yields that for $\mu$-a.e.
regular point $x$ we have $\bt (\wb ^-_\bs (x)) = \wb ^-_\bs (\bt x)  \mod 0$.
Thus we obtain the invariance of $\wb$ under the whole action $\b$.

Let $x$ and $y \in \w^-_\as (x)$ be in a Pesin set $\Rel$. We denote $x_0=x$ and
$$ x_n = \b(\s ,x_{n-1})= \b (n\s,x)  = \a (\s_n) x , \; \text{ where }\s_n =n\s + g(x,n\s)\nor.
 \eqno (1) $$
Since points $y$ and $x_n$, $n\ge1$, are in the same orbit-stable leaf
$\O \oplus \w^-_\as (x)$ we can define $y_n$ to be the intersection of the orbit
of $y$ with $\w^-_\as (x_n)$. Since all points $y_n$, $n\ge1$, are on the orbit of
$y$, we can represent $y_{n+1}$ as $\b (\s + \t_{n}, y_{n})$ for some $\t _n \in \rk$
and write
$$
y_{n} =\b (\s + \t_{n-1}, y_{n-1}) = ... = \b (n\s + (\t_0+...+\t_{n-1}), y) = \a (\s_n) y. \eqno (2)
$$
The last equality follows from invariance of $\w^-_\as$ under $\a$ which gives
that $\a (\s_n) y$ is on $\w^-_\as$ leaf of $x_n = \a (\s_n) x$ and thus coincides
with $y_{n}$ by definition. Recall that by Proposition \ref{adjust} the function $g$
satisfies
$$
|g(z,\t)| \le 2\e \| \t \| \eqno (3)
$$
for any regular point $z$ and $\t \in \rk$. Hence the sequence $\s_n =n\s + g(x,n\s)\nor$ remains in a narrow cone around the direction of $\s$. We conclude that diffeomorphisms
$\a (\s_n)$ contract the stable manifold $\w^-_\as (x)$ exponentially and thus
$$
\dist (x_n,y_n)= \dist (\a (\s_n) x, \a (\s_n) y) \le K_1 e^{- n \chi} \, \dist (x,y) \eqno (4)
$$
for some $\chi>0$ which can be chosen close to the slowest contraction rate for $\as$.

Next we will show that the series $\t = \sum_{i=0}^\infty \t_i$ converges exponentially
so that according to (2) we have $\dist (y_n, \b (n\s + \t, y)) \to 0$ exponentially.
Combining this with (4) we obtain that for $\tilde y =\b (\t,y)$
$$
\dist (\b (n\s,x),\b (n\s, \tilde y)) = \dist (x_n, \b (n\s + \t, y) ) \to 0  \eqno (5)
$$
exponentially and thus $\tilde y$ belongs to the stable ``leaf" $\wb ^-_\bs (x)$.

To show that the series $\t = \sum_{i=0}^\infty \t_i$ converges we estimate
$\t_n$ as follows. Similarly to the last equalities in (1) and (2) we can write
$$x_{n+1} = \a (\s + g(x_n,\s)\nor) x_n  \quad  \text{and } \quad
y_{n+1}=\a (\s + g(x_n,\s)\nor) y_n.  \eqno (6)
$$
Denoting $\t_n'=(g(x_n,\s) - g(y_n,\s)) \nor$
we obtain using (2) and (6) that
$$
\b (\t_n) \b(\s, y_n)=\b (\s + \t_n, y_n) = y_{n+1} = \a (\s + g(x_n,\s)\nor) y_n =
$$
$$
 \a (\t_n'+ \s + g(y_n,\s)\nor) y_n = \a (\t_n') \a(\s + g(y_n,\s)\nor) y_n = \a (\t_n') \b(\s, y_n)
$$
This shows that $\t_n$ is uniquely determined by the following equations
$$
\a (\t_n') z_n=\b (\t_n) z_n \quad \text{or} \qquad
\t_n + g(z_n,\t_n)\nor = \t_n' = (g(x_n,\s) - g(y_n,\s)) \nor,
$$
where $z_n =\b(\s, y_n)$. Using (3) we conclude that $\t_n$ is a vector parallel
to $\nor$ whose length satisfies
$$
c \, t_n' \le \| \t_n \|_{\rk} \le C \, t_n', \quad \text {where} \qquad t_n'=|g(x_n,\s) - g(y_n,\s)|.
$$
Thus we need to investigate the convergence of the series
$$
t'  = \sum_{n=0}^\infty t_n' = \sum_{n=0}^\infty \left|g(x_n,\s) - g(y_n,\s)\right|.
$$
By Proposition \ref{time change Holder} the function $g$ is H\"older continuous
with exponent $\gamma$ and constant depending on the Pesin set.
For $x$ and $y$ in the Pesin set $\Rel$, $x_n$ and $y_n $ are in another
Pesin set $\Re_\e ^{l'}$ for which the H\"older constant deteriorates from that of
$\Rel$ by a factor at most $\exp (p \e \| \s_n \|) \le \exp (2p\e n)$ (see Remark after
Proposition \ref{Lyapunov metric Holder}). Replacing $\s$
by its multiple if necessary we may assume without loss of generality that $\chi > 2p\e$.
Then (4) implies that the series converges exponentially and its sum satisfies
$$t' \le K_2 \dist (x,y)^\gamma.$$
This completes the proof of (5) and shows that $\tilde y \in \wb ^-_\bs (x)$, where
$\tilde y =\b (\t,y)$ with $\t=t\nor$ and $\| \t \| \le C \, t' \le K_3 \dist (x,y)^\gamma$.
Hence $\tilde y$ can be represented as
$$
\tilde y = \a (\varphi_x (y) \nor) y,  \qquad \text{where} \quad
|\varphi_x (y)| \le K_4 \dist (x,y)^\gamma.
$$
We conclude that $\wb ^-_\bs (x)$ is of the form stated in the proposition.
The H\"older continuity of $\varphi_x $ with H\"older exponent $\gamma$
can be obtained similarly to the H\"older  estimate for  $\varphi_x $  in the
previous equation.  The constant $K_4$ depends on the Pesin set $\Rel$.
\qed\medskip

The corresponding unstable ``foliation" $\wb ^+_\bs$ can be obtained as
$\wb ^-_\b (-\s)$. Since the foliation $\w$ corresponding to the Lyapunov
distribution $E$ is an intersection of stable foliations we obtain the following
corollary.

\begin{corollary}  \label{new manifolds}
For the Lyapunov foliation $\w$ corresponding to the distribution $E$ of the
original action $\a$ there exists "foliation" $\wb$ invariant under
the new action $\b$ of the form described in Proposition \ref{new stable}.
\end{corollary}

The foliation $\wb$ will be referred to as the Lyapunov foliation of $\b$
corresponding to the Lyapunov distribution $E$.


\subsection{Recurrence argument for the time change}\label{ergodicitysingular}
Recall that an element $\t \in \rk$ is  generic singular if it belongs to exactly one
Lyapunov hyperplane. We consider a generic singular element $\t$ in the Lyapunov
hyperplane $L$. Our goal is to show that for a typical point $x$ the limit points of
the orbit $\b (n \t)$, $n \in \N$, contain the support of the conditional measure of
$\nu$ on the leaf $\wb (x)$. More precisely, we prove the following lemma which
is an adaptation of an argument from \cite{KS3} for the current setting.
\medskip

We say that partition $\xi _1$ is {\em coarser} than $\xi _2$ (or that $\xi _2$ {\em refines}
$\xi _1$) and write $\xi _1 < \xi _2$ if $\xi _2 (x) \subset \xi _1 (x)$ for a.e. $x$.

\begin{proposition}  \label{transitivity on leaves}
For any generic singular element $\t \in L$ the partition $\xi_{\bt}$ into ergodic
components of $\nu$ with respect to $\bt$ is coarser than the measurable
hull $\xi(\wb)$ of the foliation $\wb$.
\end{proposition}

\proof
For a generic singular element $\t$ in $L$, $\chi$ is the only non-trivial Lyapunov
exponent that vanishes on $\t$.  As in Section \ref{outline tech} we take a regular
element $\s$ close to $\t$ for which $\chi (\t)>0$ and all other non-trivial
exponents have the same signs as for $\t$. Then $E_\as^-=E_\at^-$
and $E_\as^+=E_\at^+\oplus E$. Consequently, for the action $\b$ we have
$\xi (\wb_\bs^-)=\xi (\wb _\bt^-)$ and $\xi (\wb_\bs^+) < \xi (\wb_\at^+)$.
Birkhoff averages with respect to $\bt$ of any continuous function are
constant on the leaves of $\wb^-_\bt$. Since such averages
generate the algebra of $\bt$--invariant functions, we conclude that the
partition $\xi_\bt$ into the ergodic components of $\bt$ is coarser than
$ \xi(\wb_\bt^-)$, the measurable hull of the foliation $\wb^-_\bt$.
The equality $\xi(\wb_\bs^+) = \xi(\wb_\bs^-)$ is proved in the next proposition,
so we conclude that
\begin{equation}\label{pipartitiontrick}
 \xi_\bt < \xi (\wb _\bt^-) = \xi (\wb_\bs^-)  =  \xi(\wb _\bs^+) <  \xi (\wb).
\end{equation} 
\qed

\begin{remark}
Equalities  in \eqref{pipartitiontrick}  represent the ``$\pi$-partition trick'' which first appeared 
in \cite{KS3} in the setting  of actions by automorphisms of a torus.  Absence of Lyapunov exponents {\em negatively} proportional to $\chi$ is necessary  for this argument to work.   
If this condition holds for all   exponents (other than the trivial one corresponding to the orbit directions)  the action is called  {\em totally non-symplectic (TNS)}. On the other hand,  presence of exponents {\em positively} proportional to $\chi$, e.g.  non-simplicity of $\chi$ itself, 
forces considering multidimensional  coarse Lyapunov foliations, corresponding to all exponents positively proportional to $\chi$.  Naturally one cannot hope any more to have the dichotomy
of atomic vs. absolutely continuous but nevertheless under additional assumptions the $\pi$-partition trick still works and allows to make conclusions about conditional measures. 
\end{remark}

The (long) remainder of this section is dedicated to the justification of the expression ``$\pi$-partition trick'' in our setting. 

\begin{proposition}  \label{pi partition}
$\xi(\wb_\bs^+) = \xi(\wb_\bs^-) = \pi(\bs)$, the $\pi$-partition of $\bs$.
\end{proposition}

\proof We will show that  $\xi(\wb_\bs^+) = \pi(\bs)$. The other equality
is obtained in the same way. We note that in the case of diffeomorphisms
this result is given by Theorem B in \cite{LY}. Also, for the case of a hyperbolic
measure this result was established earlier in \cite[Theorem 4.6]{L}. In our case,
although some zero Lyapunov exponents appear, they correspond to the orbit
direction, so that the central direction will be easier to control than in \cite{LY}.
We will follow \cite{L} and \cite{LY}, so let us give a sketch of the proof in
their case.


\begin{subsubsection}{Sketch of the proof in \cite{L, LY}}
In what follows, $f$ will be a diffeomorphism preserving a measure
$\mu$ with unstable foliation $\w^+$ and local unstable manifolds
$W^+$. The idea is to use the following criterium due to Rokhlin
\cite{Ro} (Theorems 12.1 and 12.3).
Given a partition $\xi$ we
denote by $\M_\xi$ the $\sigma$-algebra generated by $\xi$.
\begin{theorem}\label{rokhlin}
Let $f$ be a measure preserving transformation and assume $\xi$ is
an increasing partition, i.e $\xi > f \xi$, satisfying:
\begin{enumerate}
\item \label{1ro} $\bigvee _{n=0}^\infty f^{-n} \xi$ is the partition into
points,
\item \label{2ro} $h(f)=h(f,\xi)<\infty$.
\end{enumerate}
Then the Pinsker $\sigma$-algebra coincides mod $0$ with the
$\sigma$-algebra $\bigcap _{n=0}^\infty \M_{f^{n}\xi}$.
\end{theorem}
For two partitions $\eta$ and $\xi$,
$$H(\eta|\xi)=-\int\log \mu_x^{\xi}\left(\eta(x)\right)d\mu(x),$$
where $\mu^{\xi}$ are the conditional measures associated to the
measurable partition $\xi$. And for an increasing partition $\xi$,
$h(f,\xi)=H(f^{-1}\xi|\xi)$. Also, $H(\eta)=H(\eta|\tau)$ where
$\tau$ is the trivial partition. We shall make use of the following
known formulas of the conditional entropy. Given partitions
$\eta,\xi,\zeta$,
\begin{enumerate}
\item $H(\zeta\vee\xi|\eta)=H(\zeta|\eta)+H(\xi|\zeta\vee\eta),$
\item If $\eta>\xi$ then for any $\zeta$, $H(\zeta|\eta)\leq
H(\zeta|\xi)$,
\item $H(\zeta|\xi)\geq H(\eta|\xi)-H(\eta|\zeta)$.
\end{enumerate}
Also we shall make use of the following Lemma whose proof is left to
the reader.
\begin{lemma}\label{absentropy}
Let $\xi$ and $\eta_n$ be measurable partitions and assume that if
$D_n=\{x:\xi(x)\subset\eta_n(x)\}$ then $\mu(D_n)\to 1$. Then for
any $\zeta$,
$$
\lim H(\zeta|\xi\vee\eta_n)=H(\zeta|\xi)
$$
\end{lemma}

We say that a partition $\xi$ is subordinated to $W^+$ if $U_x
\subset \xi (x) \subset W^+(x)$ for a.e. $x$, where $U_x$ is some
open neighborhood of $x$ in $W^+(x)$.
In \cite{L, LY}  first
 a partition $\xi$  is constructed as follows,
\begin{lemma}\cite[Lemma 3.1.1.]{LY}\label{311ly}
There exists a measurable partition $\xi$ with the following
properties,
\begin{enumerate}
\item $\xi$ is an increasing partition subordinated to $W^u$,
\item $\bigvee _{n=0}^\infty f^{-n} \xi$ is the partition into
points,
\item $\bigcap _{n=0}^\infty \M_{f^{n}\xi}=\M_{\xi(\w^u)}$.
\end{enumerate}
\end{lemma}
Partitions of this type were used by Sinai \cite{S} 
to study
uniformly hyperbolic systems and were built in the general context
in \cite[Proposition 3.1.]{LS}, (see also \cite[Proposition
3.1.]{L}). 
Then it is proven that hypothesis (\ref{2ro})
of Theorem \ref{rokhlin} is satisfied by any such partition.

In \cite{L} and \cite{LY} the proof that hypothesis (\ref{2ro}) of
Theorem \ref{rokhlin} is satisfied by any of these partitions is in
various steps. On one hand, the following lemma is proven.
\begin{lemma}\cite[Lemma 3.1.2.]{LY} \label{312ly}
For any two partitions $\xi_1$ and $\xi_2$ built in Lemma
\ref{311ly}, $h(f,\xi_1)=h(f,\xi_2)$.
\end{lemma}

On the other hand, it is proven that the entropy of the partitions
$\xi$ built in Lemma \ref{311ly} approach $h(f)$. To this end, it is
built a countable partition $\P$, with finite entropy, i.e.,
$H(\P)<\infty$ and such that $h(f,\P)$ is close to $h(f)$. Then
$h(f,\P)$ and $h(f,\xi)$ are compared where $\xi$ is a partition
built in Lemma \ref{311ly}. It is in comparing this two entropies
where the proof in \cite{L} and \cite{LY} differ. In \cite{L} the
comparison follows from the properties of $\xi$ and that
$\P^+:=\bigvee_{n=0}^{\infty}f^n\P$ refines $\xi$ (this is done in
the proof of \cite[Proposition 4.5.]{L}) while in \cite{LY} more
work is needed because $\P^+$ does not a priori refines $\xi$ due to
the presence of zero exponents. However, in our case, although some
zero exponents appear, they correspond to the orbit direction, so
that $\P^+$ will essentially refine $\xi$ because of the properties
of the partition $\P$. That is why we are somehow closer to the
proof in \cite{L}. Finally, the proof ends because we can take the
partition $\P$ with entropy as close to $h(f)$ as wanted.
\end{subsubsection}


\begin{subsubsection}{Proof in our case}
Let us go now to the proof in our case. We will follow the
above sketch, but now $f:=\bs$ will not be a diffeomorphism, so we
need to take some care. Let $\wb := \wb_\bs^+$ and
$\wbl:=\wbl_\bs^+$ be the global and local unstable "manifolds"
built in Proposition \ref{new stable}. As the global "manifold"
$\wb(x)$ is a graph over $\w_\as^+(x)$, the local "manifold"
$\wbl(x)$ is the restriction of this graph to $\wl_\as^+(x)$.

As the main contraction/expansion properties of $\wb$ comes from the
contraction/expansion properties of $\w_\as^+(x)$, we will mostly
measure the distances between points in $\wb$ projecting them into
$\w_\as^+(x)$. Thus, let us define $\pi_x:\wb(x)\to\w_\as^+(x)$ the
projection and observe that the restriction of $\pi_x$ to $\wbl(x)$
is Lipschitz continuous, with Lipschitz constant depending only on
the Pesin set $x$ belongs to. Observe also that the inverse of
$\pi_x$ is not Lipschitz. Let us define for $z,y\in\w(x)$, $\tilde
d_x(y,z)=d(\pi_x(y),\pi_x(z))$.

Let us begin now with a useful Lemma,
\begin{lemma}\label{generates}
If $\eta$ is an increasing partition and $\eta(x)\subset\wbl(x)$ for
a.e. $x$ then the sequence of partitions $\{f^{-n}\eta\}$ generates,
i.e. $\bigvee _{n=0}^\infty f^{-n} \eta$ is the partition into
points.
\end{lemma}
\proof Let us see that for a.e. $x$, if $y\in (f^{-n}\eta)(x)$ for
every $n$ then $x=y$. We have that a.e. $x$ belongs infinitely many
times to some fixed Pesin set, say $R$. Take $n_i$ the sequence of
integers such that $f^{n_i}(x)\in R$. So we have that $y\in
f^{-n_i}(\eta(f^{n_i}(x)))\subset f^{-n_i}(\wbl(f^{n_i}(x)))$ for
every $i$. Since $f^{n_i}(x)$ is in $R$, we have that the projected
diameter of $\wbl(f^{n_i}(x))$ is uniformly bounded and that the
projected diameter of $f^{-n_i}(\wbl(f^{n_i}(x)))$ tends to $0$.
Hence, since $\pi_x(y)\in\pi_x(f^{-n_i}(\wbl(f^{n_i}(x))))$, the
distance between the projection of $y$ into $\wl^+_\as(x)$ with $x$
is $0$, this means that in fact $y$ is in the orbit of $x$. But this
is only possible if $x=y$ since $y\in\wbl(x)$.
\qed\medskip

Recall that by Proposition \ref{new stable} we have that $ \wb (x)
=\{ \a (\varphi_x (y) \nor) y : y \in \w^+_\as(x) \}$. Following the
above philosophy, let us say that a partition $\xi$ is subordinated
to $\wbl$ if $\xi (x) \subset \wbl(x)$ for a.e. $x$, and there is an
open neighborhood $U_x \subset\wl_\as^+(x)$ such that $\{\a
(\varphi_x (y) \nor) y : y \in U_x \}\subset \xi(x)$.

Here again, we will use the criterium in Theorem \ref{rokhlin} to
prove Proposition \ref{pi partition}, that is to prove that the
Pinsker $\sigma$-algebra coincides with the $\sigma$-algebra
generated by the "foliation" $\wb_\bs^+$.
So that let us build partitions like the ones in Lemma
\ref{311ly}.
\begin{lemma}\label{311lybis}
There exists a measurable partition $\xi$ with the following
properties:
\begin{enumerate}
\item $\xi$ is an increasing partition subordinated to $\wbl$,
\item $\bigvee _{n=0}^\infty f^{-n} \xi$ is the partition into
points,
\item $\bigcap _{n=0}^\infty \M_{f^{n}\xi}=\M_{\xi(\wb)}$.
\end{enumerate}

\end{lemma}

\proof Let us take $\hat \xi$, the measurable partition built in
Lemma \ref{311ly} for $\as$, and define the partition $\xi$ as
the graph over $\hat\xi(x)$:
$$
\xi (x) =\{ \a (\varphi_x (y) \nor) y : y \in \hat\xi(x) \}.
$$
Let us see that this is a partition that satisfies the three
properties. Property (1) follows by definition and because $\hat\xi$
satisfies also property (1). Property (2) follows from Lemma
\ref{generates}. Observe that property (3) is the same as proving
that $\bigwedge_{n=0}^{\infty}f^n\xi=\xi(\wb)$. Notice that
$\bigwedge_{n=0}^{\infty}f^n\xi$ is the graph over
$\bigwedge_{n=0}^{\infty}f^n\hat\xi$ which equals $\xi(\w^+_\as)$ by
property (3) of Lemma \ref{311ly}. So, since $\xi(\wb)$ is the graph
over $\xi(\w^+_\as)$ we get property (3).
\qed\medskip


The next step is to prove the analog of Lemma \ref{312ly}, in
fact we prove a more general result. We follow the proof
in \cite[Lemma 3.1.1.]{LY}.
\begin{lemma}\label{312lygen}
If $\xi$ is a partition built in Lemma \ref{311lybis}, and $\zeta$
is an inceasing partition such that $\zeta(x)\subset\wbl(x)$, then
$h(f,\zeta\vee\xi)=h(f,\zeta)$.
\end{lemma}
\proof For $n\geq 1$ we have
\begin{eqnarray*}
h(f,\zeta\vee\xi)&=&h(f,\zeta\vee f^n\xi)=H(\zeta\vee
f^n\xi|f\zeta\vee f^{n+1}\xi)\\
&=&H(\zeta|f\zeta\vee f^{n+1}\xi)+H(\xi|f\xi\vee f^{-n}\zeta)
\end{eqnarray*}
As $n\to\infty$, the second term goes to $0$ since \{$f^{-n}\zeta\}$
generates by Lemma \ref{generates}. So we want to show that
$H(\zeta|f\zeta\vee f^{n+1}\xi)\to H(\zeta|f\zeta)$. To this end we
shall make use of Lemma \ref{absentropy}. So let $D_n=\{x:
(f\zeta)(x)\subset (f^{n+1}\xi)(x)\}$. Since $\zeta(x)\subset
\wbl(x)$, and the projected diameter of $\wbl(x)$ into
$\wl_\as^+(x)$ is finite a.e., we have that the projected diameter
of $(f^{-n}\zeta)(x)$ goes to $0$. Hence, since $\xi(x)$ contains a
graph over an open neighborhood of $x$ in $\wl_\as^+(x)$ we have
that $(f^{-n}\zeta)(x)\subset \xi(x)$ if $n$ is big enough and hence
$\mu(D_n)\to 1$. Now the lemma follows from Lemma \ref{absentropy}
and the fact that $h(f,\zeta)=H(\zeta|f\zeta)$ since $\zeta$ is an
increasing partition. \qed

\begin{corollary}\label{312lybis}
For any two partitions $\xi_1$ and $\xi_2$ built in Lemma
\ref{311lybis}, $h(f,\xi_1)=h(f,\xi_2)$.
\end{corollary}

It remains to show that the entropy of a partition built in Lemma
\ref{311lybis} equals the entropy of $f$.
We shall build a countable partition $\P$ with finite
entropy to compare $h(f,\P)$ with $h(f,\xi)$ as in the sketch. To
this end we shall use the following Lemma due to Ma\~n\'e,
\begin{lemma}\cite[Lemma 2]{M}\label{2ma}
If $\mu$ is a probability measure and $0<\psi<1$ is such that
$\log\psi$ is $\mu$ integrable, then there exists a countable
partition $\P$ with entropy $H(\P)<+\infty$ such that $\P(x)\subset
B(x,\psi(x))$ for a.e. $x$.
\end{lemma}

So let us construct a suitable function $\psi$. For a set $A\subset
M$ let us define $\O_\e A=\{\at(a): a\in A\,;\,\|\t\|<\e\}$. Let us
fix a Pesin set $R_1$ of positive measure and take $R_0$ another
Pesin set such that $\O_\e R_1\subset R_0$. Arguing similar to Lemma
2.4.2. of \cite{LY} let us define a measurable function
$\psi:S\to\R^+$ by
\begin{displaymath}
\psi(x)=\left\{\begin{array}{lll}
\d & \mbox{if} & x\notin R_0 \\
\d l_0^{-1}e^{-\la r(x)} & \mbox{if} & x\in R_0 \\
\end{array}\right.
\end{displaymath}
Where $r(x)$ is the smallest positive integer $k>0$ such that
$f^k(x)\in R_0$, $\lambda$ and $l_0=l_{R_0}$ are the constants in
Lemma \ref{lemma242} below and $\d$ is such that if $x,y\in R_0$ and
$\dist(x,y)<\d$ then $\O_\e W^+_\as (x)\cap \O_\e W^-_\as
(y)\neq\emptyset$ and vice versa interchanging $x$ and $y$ for some
$\e>0$ small that depends on the Pesin set (such $\d$ and $\e$
exists by transversality and uniformity over Pesin sets). We will
require other properties for that $\delta$ later (see Lemma
\ref{lemma242b}). Since $\int_{R_0}rd\mu=1$, we get that $\log\psi$
is integrable. We may assume also, by an appropriate choice of
$R_0$, that $\inf_{n\geq 0} \psi(f^{-n}(x))=0$ for a.e $x$.

Hence, by Lemma \ref{2ma}, there is a partition $\tilde \P$ such
that $H(\tilde\P)<\infty$ and $\tilde\P(x)\subset B(x,\psi(x))$ for
a.e. $x$. Take $\hat R_1\subset R_1$ such that if $x, y\in\hat R_1$
and $\dist (x,y)<\d$ then there is a point $z\in R_1\cap \O_\e
W^+_\as (x)\cap \O_\e W^-_\as (y)$. If $R_1$ is taken of big enough
measure, then there is such a set $\hat R_1$ of positive measure.
Let us define $\P=\tilde\P\vee\{\hat R_1,S\setminus \hat R_1\}$ and
recall that $\P^+=\bigvee_{n=0}^{\infty}f^n\P$.

\begin{lemma}\label{lemma242b}
For some $\delta>0$ we have that $\P^+(x)\subset \wbl (x)$, $x$ a.e.
\end{lemma}
Before the proof of this Lemma, let us begin with a property of the
$\wbl$ ``manifolds" that, apart from invariance and uniformity over
Pesin sets, simply reflects the Lipschitz property of the original
map $\as$.
\begin{lemma}\label{lemma242} There is $\la>0$ that
depends on $\s$, $l=l_R>0$ that depends on the Pesin set $R$ such
that for $n>0$, and $0<\d\leq l^{-1}$, if $z\in\wbl(x)$, $x\in R$
are such that $\tilde d_x(x,z)<\d e^{-n\la}$ then $$\tilde
d_{f^n(x)}(f^n(x),f^n(z))<\d$$ and $f^n(z)\in\wbl(f^n(x))$.
\end{lemma}

Let us go now into the proof of Lemma \ref{lemma242b}.

\proof{of Lemma \ref{lemma242b}.} Let us see first that
$\P^+(x)\subset\O_\e W^+_\as (x)$, $x$ a.e. Let $y\in \P^+(x)$. Take
the sequence of negative integers $-n_i<-n_{i-1}$ such that
$x_{n_i}=f^{-n_i}(x)\in\hat R_1$. Since
$y\in\P^+(x)=\bigvee_{n=0}^{\infty}f^n\P$ and $\P=\tilde\P\vee\{\hat
R_1,S\setminus \hat R_1\}$ we have that $y_{n_i}=f^{-n_i}(y)\in\hat
R_1$ and hence, since $\dist(x_{n_i},y_{n_i})<\d$, we get that
$$R_1\cap \O_\e W^+_\as (x_{n_i})\cap \O_\e W^-_\as
(y_{n_i})\neq\emptyset.$$ So, the whole piece of orbit $\O_\e
W^+_\as (x_{n_i})\cap \O_\e W^-_\as (y_{n_i})$ is in $R_0$. Call
$z^1_{n_i}=\wbl^+_\bs (x_{n_i})\cap \O_\e W^-_\as (y_{n_i})$ and
$z^2_{n_i}=\O_\e W^+_\as (x_{n_i})\cap \wbl^-_\bs (y_{n_i})$.

We claim that $f^{n_i-n_{i-1}}(z^1_{n_i})=z^1_{n_{i-1}}$ for every
$i$ and hence $f^{n_i-n_j}(z^1_{n_i})=z^1_{n_j}$. The same happens
for the sequence $z^2_{n_i}$. Let us proof the claim.

Since $z^1_{n_i}= \wbl^+_\bs (x_{n_i})\cap\O_\e W^-_\as (y_{n_i})$
we have that
$$f^{n_i-n_{i-1}}(z^1_{n_i})=f^{n_i-n_{i-1}}(\wbl^+_\bs
(x_{n_i}))\cap f^{n_i-n_{i-1}}(\O_\e W^-_\as (y_{n_i})).$$ Hence to
prove the claim it is enough to show that
$f^{n_i-n_{i-1}}(z^1_{n_i})\in \wbl^+_\bs (x_{n_{i-1}})$ and to this
end we shall use Lemma \ref{lemma242}. Take the sequence of positive
integers $k_j, j=0,\dots l$ such that $f^{k_j}(x_{n_i})$ enters in
$R_0$, $k_0=0$, $k_l=n_i-n_{i-1}$. By definition we have that
$k_j-k_{j-1}=r(f^{k_{j-1}}(x_{n_i}))=r_{j-1}$. Let us see that
$$f^{r_{j-1}}(f^{k_{j-1}}(z^1_{n_i}))\in\wbl^+_\bs
(f^{r_{j-1}}(f^{k_{j-1}}(x_{n_i}))).$$ By Lemma \ref{lemma242} it is
enough to see that
$$\tilde d_{f^{k_{j-1}}(x_{n_i})}(f^{k_{j-1}}(x_{n_i}),f^{k_{j-1}}(z^1_{n_i}))<l_0^{-1} e^{-r_{j-1}\la}.$$
Let us assume by induction that
$$f^{k_{j-1}}(z^1_{n_i})=\wbl^+_\bs
(f^{k_{j-1}}(x_{n_i}))\cap\O_\e W^-_\as (f^{k_{j-1}}(y_{n_i})).$$
Since $f^{k_{j-1}}(x_{n_i})\in R_0$ we know that
$$\dist(f^{k_{j-1}}(x_{n_i}),f^{k_{j-1}}(y_{n_i}))<\d e^{-r_{j-1}\la}.$$
Now, by the uniformity of the invariant stable and unstable
manifolds for points in a given Pesin set and by the uniform
transversality of the invariant distribution, there is a constant
$C_0$ that depends on the Pesin set such that
$$\tilde d_{f^{k_{j-1}}(x_{n_i})}(f^{k_{j-1}}(x_{n_i}),f^{k_{j-1}}(z^1_{n_i}))\leq C_0\dist(f^{k_{j-1}}(x_{n_i}),f^{k_{j-1}}(y_{n_i})).$$
So that taking $\d$
small enough we get that
$$f^{r_{j-1}}(f^{k_{j-1}}(z^1_{n_i}))\in \wbl^+_\bs
(f^{r_{j-1}}(f^{k_{j-1}}(x_{n_i}))),$$ hence
$$f^{k_j}(z^1_{n_i})=\wbl^+_\bs
(f^{k_j}(x_{n_i}))\cap\O_\e W^-_\as (f^{k_j}(y_{n_i})).$$ So the
claim is proved for $z^1_{n_i}$.

For the case of $z^2_{n_i}$ observe that $z^2_{n_i}=\O_\e
(z^1_{n_i})\cap \wbl^-_\bs (y_{n_i})$ for every $i$. On the other
hand, $\wbl^-_\bs$ is $f$-invariant and since the derivative of
$f^n$ restricted to any orbit $\O$ is uniformly bounded from below
and from above we get that
\begin{eqnarray*}
f^{n_i-n_{i-1}}(z^2_{n_{i-1}})&=&f^{n_i-n_{i-1}}(\O_\e
(z^1_{n_{i-1}}))\cap f^{n_i-n_{i-1}}(\wbl^-_\bs
(y_{n_{i-1}}))\\
&\subset& \O_{C\e}(f^{n_i-n_{i-1}}(z^1_{n_{i-1}}))\cap \wbl^-_\bs
(f^{n_i-n_{i-1}}(y_{n_{i-1}}))\\
&=&\O_{C\e}(z^1_{n_i})\cap \wbl^-_\bs (y_{n_i})
\end{eqnarray*}
for some fixed constant $C$. So, if $\e$ is small enough we get that
the last term equals $z^2_{n_i}$.

So finally, since $f^{n_i-n_0}(z^2_{n_i})=z^2_{n_0}$, we get that
$$\dist(z^2_{n_0}, y_{n_0})=
\dist(f^{n_i-n_0}(z^2_{n_i}),f^{n_i-n_0}(y_{n_i}))$$ and hence,
since the right hand side tends to zero because $z^2_{n_i}\in
\wbl^-_{\bs}(y_{n_i})$ we get that $z^2_{n_0}=y_{n_0}$. Hence,
$f^{-n_0}(y)\in\O_\e W^+_\as (f^{-n_0}(x))$ and since $n_0\leq
r(f^{-n_0}(x))$ we get that $\P^+(x)\subset\O_\e W^+_\as (x)$ by
using Lemma \ref{lemma242} and the fact that the derivative of $f^n$
restricted to any orbit $\O$ is uniformly bounded from below and
from above.

So we get that $\P^+(x)\subset\O_\e W^+_\as (x)$, $x$ a.e. Let us
see now that in fact $\P^+(x)\subset \wbl^+_\bs (x)$, $x$ a.e.
We will use the same notations as above. Take $y\in \P^+(x)$, we
already get that $y_{n_0}=z^2_{n_0}$ and hence that it is in the
$\e$-orbit of $z^1_{n_0}$, where $z^1_{n_0}=\wbl^+_\bs (x_{n_0})\cap
\O_\e W^-_\as (y_{n_0})$. Let us see that $z^1_{n_0}=y_{n_0}$. In
fact, $\dist(f^{-n}(z^1_{n_0}),f^{-n}(x_{n_0}))\to 0$ since
$z^1_{n_0}\in\wbl^+_\bs (x_{n_0})$. On the other hand, since
$\inf_{n\geq 0}\psi(f^{-n}(x))=0$, we get that $\liminf
\dist(f^{-n}(y_{n_0}),f^{-n}(x_{n_0}))=0$. So $$\liminf
\dist(f^{-n}(y_{n_0}),f^{-n}(z^1_{n_0}))=0.$$ But the derivative of
$f^n$ restricted to any orbit $\O$ is uniformly bounded from below
and from above, that is: $\|D_zf^n|_{T\O }\|\leq C$ for every
$n\in\Z$. So we have that $$C^{-1}\dist(y_{n_0},z^1_{n_0})\leq
\dist(f^{-n}(y_{n_0}),f^{-n}(z^1_{n_0}))$$ for every $n$ and hence
$\dist(y_{n_0},z^1_{n_0})=0$. Hence, $f^{-n_0}(y)\in \wbl^+_\bs
(f^{-n_0}(x))$ and since $n_0\leq r(f^{-n_0}(x))$ we get the Lemma
by using Lemma \ref{lemma242}.
\qed\medskip

So we can now begin the comparison of the entropies $h(f,\P)$ and
$h(f,\xi)$ for the partition $\P$ built just before Lemma
\ref{lemma242b} and the partition $\xi$ built in Lemma
\ref{311lybis}. But first let us state the following corollary of
Lemma \ref{312lygen}.

\begin{corollary}\label{corp0}
Let $\P$ be the partition in Lemma \ref{lemma242b} and $\Qe$ be any
finite entropy partition. Then, for $\P_0=\P\vee\Qe$ and $\xi$ a
partition built in Lemma \ref{311lybis} we have that
$h(f,\P_0)=h(f,\P_0^+)=h(f,\xi\vee\P_0^+)$.
\end{corollary}


\proof The result follows since $\P_0^+(x)\subset\wbl(x)$. \qed\medskip

Finally we get:
\begin{lemma}\label{321ly}
$h(f,\xi\vee\P_0^+)=h(f,\xi)$.
\end{lemma}
\proof As in the argument in the proof of Lemma \ref{312lygen}, see
also \cite[Lemma 3.2.1.]{LY}, we have,
\begin{eqnarray*}
h(f,\xi\vee\P_0^+)&=&h(f,\xi\vee f^n\P_0^+)\\
&=&H(\xi|f\xi\vee f^{n+1}\P_0^+)+H(\P_0^+|f^{-n}\xi\vee f\P_0^+)
\end{eqnarray*}
where the first term is $\leq H(\xi|f\xi)$ and the second term goes
to $0$ since $\{f^{-n}\xi\}$ generates. Hence
$h(f,\xi\vee\P_0^+)\leq H(\xi|f\xi)=h(f,\xi)$. Finally, since
$H(\P_0)<\infty$ we have that $h(f,\xi\vee\P_0^+)\geq h(f,\xi)$ and
hence we are done. \qed\medskip

Finally, combining the above Lemma with Corollary \ref{corp0} we get
that $h(f,\P_0)=h(f,\xi)$. So that taking finite partitions $\Qe_n$
such that $h(f,\Qe_n)\to h(f)$ we get that
$$h(f,\xi)=h(\P\vee\Qe_n)\geq h(f,\Qe_n)$$
and hence $h(f,\xi)\geq h(f)$. The other enequality follows since
$\xi$ is a measurable partition.

%
%
%
%
%

\end{subsubsection}


\section{Conclusion of the proof of Theorem \ref{tech}.}\label{sectionconclusion}

We will use the properties of the time change and the transitivity property
of the action $\b$ to produce  elements of the action $\a$
with recurrence and uniformly bounded derivatives along $\w$.

We denote by $\mw _x$ the conditional measure of $\mu$ on $\w (x)$ and
by $B^ \w _r (x)$ the ball in $\w (x)$ of radius $r$ with respect to the induced
smooth metric.

\begin{lemma}  \label{transitivity for old}
For any Pesin set $\Rel$ there exist positive constants $K$ and $l'$ so that
for $\mu$- a.e. $x \in \Rel $ and for $\mw _x$- a.e. $y \in \Rel \cap B^ \w _r (x)$
there exists a sequence of elements $\t _j \in \rk$ with
\begin{enumerate}
\item $x_j = \a (\t _j)x  \in \Re _\e ^{l'}$,
\item $x_j \to y$,
\item $K^{-1} \le \| \de _x \a (\t _j) \| \le K$.
\end{enumerate}
\end{lemma}

\proof
Consider typical points $x \in \Rel $ and $y \in \Rel \cap B^ \w _r (x)$
and let $\tilde y =\a (\varphi_x (y) \nor) y $ be the point on $\wb (x)$
corresponding to $y$. We denote $\s = \varphi_x (y) \nor$ and observe that
$\tilde y \in \Re _\e ^{l''}$ with  $l'' = l \exp (\|\s \|)$. By Proposition \ref{new stable}
the function $\varphi_x (y)$ is H\"older on $\Rel \cap B^ \w _r (x)$ and hence
$\s$ is uniformly bounded, so that the constant $l''$ can be chosen the same
for all $x$ and $y$ in the lemma.

As we show below, $x$ and  $\tilde y$ are are also typical points with respect 
to the invariant measure $\nu$ for $\b$. Then by
Proposition \ref{transitivity on leaves} there exists a sequence $n_j \to \infty$
such that $\tilde x_j = \b(n_j \t,x) =\a (\g(x, n_j \t))x \to \tilde y$.
Since both $x$ and $\tilde y$ are in $\Re _\e ^{l''}$ then the iterates $\tilde x_j$
can also be taken in this set. Denoting $\t _j = \g(x, n_j \t) - \s$ we conclude
that $x_j = \a (\t _j) x \to y$. Again, all points $x_j$ are in a Pesin set $\Re _\e ^{l'}$
with $l'$ the same for all $x$ and $y$ in the lemma. Thus the sequence $\t _j$
satisfies (1) and (2).
To obtain (3) we note that by the definition of the time change
$$\| \de _x \a (\g(x, n_j \t) \|_{\e} =1.$$ Then the estimates in (3) follow from the
uniform boundedness of the correction $\s$ for all $x$ and $y$ in the lemma
and from the uniform estimates \eqref{Lyapunov>metric},\eqref{Lyapunov<metric} for
the ratio of the Lyapunov and smooth norms on the Pesin set $\Re _\e ^{l'}$.

We will now show that $x$ and  $\tilde y$ are $\nu$-typical. Since measures 
$\mu$ and $\nu$ are equivalent we may assume that $x$ is $\nu$-typical. 
It remains to prove that $\tilde y$ is typical for the conditional measure 
$\nu _x ^{\wb}$ of $\nu$ on $\wb$. For this we need to show that the
holonomy along $\on$ between the leaves $\w (x)$ and $\wb (x)$ 
is absolutely continuous with respect to the measures $\mw _x$ and $\nu _x ^{\wb}$.
We consider the foliation $W=(\on \oplus\w)=(\on \oplus\wb)$.
We note that the conditional measures $\mu_x$ and $\nu_x$ of $\mu$ and $\nu$ 
on $W(x)$ are also equivalent. Since $\on$ is the orbit foliation of the 
one-parameter subgroup $\{t \nor \}$ for both $\a$ and $\b$, the
the conditional measures $\mu_x$ and $\nu_x$ are locally equivalent to
the product of $\mw _x$ and $\nu _x ^{\wb}$ with the conditional measures
on $\on(x)$ for $\mu$ and $\nu$ respectively. The latter measures are 
equivalent to Lebesgue on $\on (x)$, for $\nu$ this follows from differentiability
of the time change $\b$ along the orbits. Since the time change $\b$, as well as 
the leaf $\wb (x)$ viewed as a graph over $\w (x)$, is also continuous on Pesin 
sets, it follows that the holonomy along $\on$ is absolutely continuous.
\qed
\medskip

We will use the notion of an affine map on a leaf of a Lyapunov foliation. These are
the maps which are affine with respect to the atlas given by affine parameters on
these leaves. The notion of affine parameters is similar to that of nonstationary
linearization. The following proposition provides $\a$-invariant affine parameters
on the leaves of any Lyapunov foliation $\w$.

\begin{proposition}       \label{affine structures} \cite[Proposition 3.1, Remark 5]{KaK3}
There exists a unique measurable family of $\Ce$ smooth $\a$-invariant affine
parameters on the leaves $\w (x)$. Moreover, they depend uniformly
continuously in $\Ce$ topology on $x$ in any Pesin set.
\end{proposition}

Now we can apply Lemma \ref{transitivity for old} to obtain the following
invariance property for the conditional measures of $\mu$ on $\w$.
We note that if the conditional measures $\mw_x$ are atomic
this invariance property degenerates into trivial.

\begin{lemma}  \label{transitive group}\cite[Lemma 3.9]{KaK3}
For $\mu$- a.e. $x \in \Rel$ and for $\mw _x$- a.e.  $y \in \Rel \cap B^ \w _r (x)$
there exists an affine map $g : \w (x) \to \w (x)$ with $g (x)=y$ which preserves the
conditional measure $\mw _x$  up to a positive scalar multiple.
\end{lemma}

This lemma is proved by finding a limit for the restrictions of maps $\a (\t_j)$ to $\w (x)$.
The proof is identical to the one of Lemma 3.9 in \cite{KaK3}. It relies only on the
conclusions of Lemmas 3.7 and 3.8 in \cite{KaK3}, which are now given by
Lemma \ref{transitivity for old}.

Assuming that the conditional measures $\mw_x$ are non-atomic for
$\mu$-a.e. $x$, the following lemma from \cite{KaK3} establishes the
absolute continuity of these conditional measures and completes the
proof of Theorem \ref{tech}. Its proof in \cite{KaK3} relies only on the
conclusion of Lemma \ref{transitive group} (Lemma 3.9 in \cite{KaK3}).

\begin{lemma}  \label{conditional abs cont}  \cite[Lemma 3.10]{KaK3}
The conditional measures $\mw_x$  are absolutely continuous for
$\mu$ - a.e. $x$. (In fact, $\mw_x$ is Haar with respect to the affine
parameter  on $\w (x)$.)
\end{lemma}


\section{Concluding remarks and some open problems}

\subsection{Further properties  of maximal rank actions}
For a $\Z^k, \,\,k\ge 2$  action $\a$ on the  torus $\T^{k+1}$  with Cartan homotopy data  there is a unique invariant measure $\mu$ which is projected to Lebesgue measure $\lambda$ by the semi-conjugacy with the corresponding linear Cartan action $\a_0$; this measure is absolutely continuous and the semi-conjugacy is bijective  and measure preserving between certain sets of full  $\mu$-measure and full Lebesgue  measure. Thus $(\a, \mu)$  and $(\alpha_0, \la)$  are isomorphic   as measure preserving actions (\cite[Corollary 2.2.]{KRH}); furthermore, the  measurable conjugacy is smooth  on almost every local (and hence global)   stable manifold for any element of $\alpha$, in particular, along the Lyapunov foliations, \cite[Proposition 2.9.]{KRH}.  This implies that  the Jacobians along those foliations are rigid, i.e. multiplicatively cohomologous to the corresponding eigenvalues of elements of $\a_0$ (\cite[Lemma 4.4]{KaK3})\footnote{In fact, we can prove in this setting rigidity of general cocycles  which are H\"older  with respect to the Lyapunov metric. Although the proof is not very difficult it uses  the semi-conjugacy and its regularity properties  very heavily and thus  it would not fit well  with the program of the present paper which aims at deriving  geometric/rigidity  properties from purely dynamical  assumptions irrespective of any model. }. Another consequence is that the metric  entropy     of $\a $ with respect to $\mu$  is the logarithm of an  algebraic integer of degree at most $k+1$.

It is natural to ask  whether in our more general setting  similar properties of the expansion coefficients and for entropy hold. In a  general setting  Jacobian along a Lyapunov foliation is called {\em rigid} if its logarithm is cohomologous  (with a measurable transfer function) to the corresponding Lyapunov exponent. Notice that our proof of the key recurrence property  is based on  rigidity of Jacobians for the special time changes constructed in Section~\ref{measurable time change} which is true essentially by definition. Notice however that different Lyapunov foliations may require different time changes.

\begin{conjecture}\label{Jacobianrigidity} Jacobians along  Lyapunov foliations for an  action $\a$ satisfying assumptions of Main Theorem are rigid.
\end{conjecture}

\begin{problem} \cite{KaOb} What are possible values of entropy for elements of an  action  satisfying assumptions of Main Theorem?
\end{problem}

The following conjecture represents a cautiously  optimistic view  which presumes existence of a certain  underlying arithmetic structure.

\begin{conjecture} The  values of  Lyapunov exponents and  hence of entropy are  logarithms of  algebraic integers. \end{conjecture}

Notice that this is true for all known  examples on a variety of manifolds as described in the introduction and in those cases the algebraic  integers  have degree at most $k+1$.

Another result of \cite{KRH} (Theorem 3.1.) establishes existence of a  set of periodic points dense in the support of the measure $\mu$ whose eigenvalues are equal to the corresponding powers of the eigenvalues of $\a_0$.   This implies that  the Lyapunov exponents of  atomic measures concentrated on the corresponding periodic orbits are equal to those of $\mu$.   Again the latter property is also true for all known examples on manifolds other than tori. Let us call such periodic points  {\em proper}.

\begin{problem}
 Under the assumptions of Main Theorem (1) is there any proper periodic point for $\a$? Are proper periodic points dense in the support of the measure $\mu$?
 \end{problem}
 
 Another circle of questions concerns  relations  between the $\Zk$ actions  satisfying assumptions of the Main Theorem (1) and $\rk$ actions satisfying assumptions of Main Theorem (2).  Any suspension of the action of the first kind  is an action of the second kind. 
 One can also make  time changes for the suspension. A trivial type of time change is given by a linear automorphism of $\rk$.  Any time change is given by an $\rk$ cocycle over the action and Lyapunov exponents are transformed according to the cocycle averages and hence assumptions  of Main Theorem (2) are preserved under a smooth time change. Thus existence of non-trivial  time changes is  closely related to the problem of cocycle rigidity. 
 \begin{problem} Is any smooth $\R$-valued cocycle over an action satisfying assumptions
 of Main Theorem (1) (or (2))  cohomologous to a constant cocycle?
 \end{problem}
 
 The answer may depend on the  regularity of cohomology. In particular, it is more likely to be positive if the cohomology in question is only measurable, rather than smooth. 
 
 \begin{problem}Are there  $\rk$ actions satisfying assumptions of Main Theorem (2) which do not appear from time changes of suspensions of $\Zk$ actions satisfying assumptions of Main Theorem (1)?
 \end{problem}
 
 Notice that  on the torus for an action with Cartan homotopy data the unique ``large'' invariant measure, ie. the measure which projects to Lebesgue measure for the linear Cartan action  under the semi-conjugacy  
 (see\cite{KRH})  changes continuously in  weak*  topology.  Thus there is not only global but also local rigidity for     such a measure.  While global rigidity is problematic in the setting of Main Theorem,
the local version is  plausible.

\begin{conjecture} Given a $C^2$ action $\a$ with an invariant measure $\mu$ satisfying assumptions of the Main Theorem (1) or (2), 
any action $\a'$ close to $\a$ in $C^2$ topology has an ergodic invariant  measure $\mu'$ satisfying the same assumptions.  One can chose $\mu'$ in such a way that  when $\a'$ converges to $\a$  in $C^2$ topology, $\mu'$ converges to  $\mu$ in weak* topology.

Furthermore, in the  $\Zk$ case Lyapunov exponents of  $\mu'$ are equal to those of $\mu$. 
\end{conjecture}

\subsection{High rank and low dimension}

As  explained in \cite[Section 4]{KRH} many examples of manifolds with actions satisfying assumptions of  Main
Theorem (1) can be obtained  by  starting from the torus and  applying two procedures:
\begin{itemize}
\item blowing up  points and glueing in copies of the projective space $\R P(k)$, and
\item cutting  pairs of holes  and attaching handles $S^k\times\mathbb D^1$.
\end{itemize}

\begin{conjecture} An action satisfying assumptions of   Main Theorem  (1) exists on any  compact manifold of dimension three or higher.
\end{conjecture}

The sphere $S^3$ looks as a good open test  case.

\begin{definition}\cite{K97} An  ergodic invariant measure of a  $\Zk$ action with non-vanishing Lyapunov exponents
is called {\em strongly hyperbolic} if   intersection of all Lyapunov hyperplanes is the origin.
\end{definition}
Obviously the rank of a strongly hyperbolic action does not exceed the dimension of the manifold.
Furthermore, any ergodic measure for a strongly hyperbolic action of $\Zk$ on a $k$-dimensional manifold is atomic and  is supported by a single closed orbit, \cite[Proposition 1.3]{K97}. Thus   the maximal rank for a strongly hyperbolic  action
on a manifold $M$ with a  non-atomic ergodic measure, in particular measure with positive entropy, is  $\dim M-1$, exactly the case considered in the present paper.

Let us consider the  lowest dimension compatible with the higher rank assumption, namely  strongly hyperbolic $\Z^2$ actions on three-dimensional manifolds.
Lyapunov hyperplanes  are lines in this case and  the general position condition is equivalent to three  Lyapunov lines being different.   In this case our theorem applies and any ergodic  invariant measure either has  zero entropy for all elements of $\Z^2$ or is absolutely continuous.

Let us consider other possible configurations of Lyapunov lines:
\begin{enumerate}
\item Two Lyapunov exponents proportional with negative proportionality coefficient; two Lyapunov lines.
\item Two Lyapunov exponents proportional with positive proportionality coefficient; two Lyapunov lines.
\item All three Lyapunov exponents proportional; one Lyapunov line.
\end{enumerate}
\medskip

(1)  First notice that such a measure cannot be absolutely continuous. For an absolutely continuous invariant measure sum of the Lyapunov exponents is identically equal to zero. 
But in this case  two exponents are zero along the   common kernel of two  proportional exponents while the third one is not zero there. 

Now we construct an example of an action with a  singular positive entropy measure of this type. 
Consider the following  action on $\T^3$: Cartesian product of the  action generated by a diffeomorphism  $f$ of $S^1$ with one contracting  fixed point $p$ with positive eigenvalue $\beta< 1$ and one expanding
fixed point,  with  the action generated  by a hyperbolic automorphism $F$ of $\T^2$ with an    eigenvalue  $\rho>1$. The
measure $\mu=\delta_p\times\lambda_{\T^2}$ is invariant under the Cartesian product
and is not absolutely continuous.
 Lyapunov exponents are $x\log\beta ,\,y\log\rho$ and $-y\log \rho$,
 the entropy is $h_{\mu}(f^mF^n)=|n|\log\rho$.
 \smallskip

(2)  There are four Weyl chambers and in one of those all three Lyapunov exponents are negative; hence by  \cite[Proposition 1.3]{K97} any ergodic invariant measure  is atomic.
Notice that this includes the case of  a multiple exponent.
 \smallskip

 (3) On the Lyapunov line all three Lyapunov exponents vanish; hence the action is not
 strongly hyperbolic.
 \bigskip

Thus we obtain the following  necessary and sufficient condition for a configuration of Lyapunov lines.
\begin{corollary}
A  strongly hyperbolic ergodic  invariant measure  for a $\Z^2$ action on a three-dimensional
compact manifold with positive entropy  for some element is absolutely continuous if and only if  Lyapunov lines for three  exponents are different.

\end{corollary}

There is an open question related to the case (1).
In our examples both Lyapunov lines are rational. One can modify this example to make the  ``single'' Lyapunov line (the kernel of a single exponent) irrational. 
It is conceivable, although  not very likely,
that  situation when the ``double'' line, i.e. the kernel of two exponents,  or both  lines, are irrational may be different.

\begin{problem} Construct an example of a smooth $\Z^2$ action on a compact 
three-dimensional manifold with a singular ergodic invariant measure with positive entropy with respect to some element of the action, such that  two Lyapunov exponents are negatively proportional and their common kernel is an irrational line. 
\end{problem}

\subsection{Low rank and high dimension} Essentially all known rigidity results for
{\em algebraic actions} (hyperbolic or partially hyperbolic), including  cocycle,
measurable, and local differentiable rigidity, assume only some sort of ``genuine
higher ($\ge 2$)  rank'' assumption, see e.g. \cite{KS, KKS, DK-KAM} and references thereof.
By contrast, global rigidity for Anosov actions on a torus \cite{RH} and nonuniform
measure rigidity on a torus \cite{KaK3,KRH}, as well as  results of this paper, deal
with maximal  rank actions. Notice, however, that global rigidity results for Anosov
actions on an arbitrary manifold satisfying stronger dynamical assumptions only
require rank $\ge 3$ \cite{KaSp} or rank $\ge 2$ \cite{KaSa3,KaSa4}

We expect that global measure rigidity results, both on the torus  and in the general setting similar to those of the present paper, hold  at  a greater generality  although  we do not see a realistic approach for the most general ``genuine higher rank''  situation even on the torus. There is  still an intermediate class  which  is compatible with the lowest admissible rank (i.e rank  two) on manifolds of arbitrary dimension. We already mentioned it in the remark in Section~\ref{ergodicitysingular}.

 \begin{definition} An ergodic  invariant measure for a $\Zk$ action is called {\em totally non-symplectic (TNS)} if  for any two Lyapunov exponents there exists an element of $\Zk$  for which  both those exponents are negative.
 \end{definition}

Equivalently all Lyapunov  exponents are non-zero and there are no  proportional exponents with negative coefficient of proportionality\footnote{The latter is exactly what happens for symplectic actions; hence the name.}.

The TNS condition is the most general one  for the ``$\pi$-partition trick'' to work.
It also greatly helps in the geometric treatment of cocycle rigidity, see e.g. \cite{KNT}.
It has a nice property that it is inherited  by a restriction of the action to a subgroup  of rank $\ge 2$ if it is in general position.  While it is possible that  our result generalizes to the TNS measures, i.e.  assuming  only existence of some elements with positive entropy,  we  prefer to be more conservative  and formulate a conjecture under   a stronger entropy assumption.

\begin{conjecture} Let $\mu$ be an ergodic invariant totally non-symplectic measure for a smooth action $\a$
of  $\Zk,\,\,k\ge 2$. Assume that every element other than identity  has positive entropy.
Then $\mu$ is absolutely continuous.
\end{conjecture}

A serious difficulty  even in TNS case may appear in the presence  of multiple exponents.
Recall that even for linear actions  multiple eigenvalues lead to  Jordan blocks so that
when the eigenvalue  has absolute value one the action is not isometric.  More generally, positively proportional eigenvalues   also lead to complications.  Thus a  more tractable case   would  be that  with simple Lyapunov exponents  and no proportional ones. In this case the Lyapunov distributions  are one--dimensional and coincide with coarse Lyapunov ones.  For the suspension every Lyapunov exponent satisfies assumptions  of Theorem~\ref{tech}. This is a nonuniform counterpart of Cartan actions in the sense of \cite{KaSp}.

\begin{conjecture} Let $\mu$ be an ergodic invariant measure for a smooth action $\a$
of  $\Zk,\,\,k\ge 2,$ such that all Lyapunov exponents are simple and  all Lyapunov hyperplanes different.  Assume that some  element of the action  has positive entropy.
Then $\mu$ is absolutely continuous.
\end{conjecture}



{\small

\end{document}